\documentclass[11pt]{article}
\usepackage[english]{babel}
\usepackage[a4paper,scale=0.75,hcentering]{geometry}
\usepackage{amsmath,amsthm,amssymb,mathtools,booktabs,enumitem}
\usepackage{xcolor}
\usepackage[hyperfootnotes=false]{hyperref}
\hypersetup{colorlinks=true,linkcolor=blue,urlcolor=red,citecolor=blue,
 pdftitle={Sharp Sobolev inequalities on the sphere: optimal powers and remainder estimates},
 pdfauthor={Yuxuan Zhou and Wenming Zou}}
\allowdisplaybreaks
\newcounter{results}[section]

\theoremstyle{plain}
\newtheorem{theorem}[results]{Theorem}
\newtheorem{lemma}[results]{Lemma}
\newtheorem{proposition}[results]{Proposition}
\newtheorem{corollary}[results]{Corollary}
\numberwithin{equation}{section}
\newcommand{\R}{\mathbb R}
\newcommand{\D}{\mathcal D}
\DeclareMathOperator{\Var}{Var}
\title{\bfseries Sharp Sobolev inequalities on the sphere: optimal powers and remainder estimates\thanks{This work is supported by National Key R\&D Program of China (Grant 2023YFA1010001) and NSFC (Grant 12571123, 12671139).\newline E-mail addresses: zhou-yx22@mails.tsinghua.edu.cn (Zhou), zou-wm@mail.tsinghua.edu.cn (Zou).}}
\author{{\bfseries Yuxuan Zhou, Wenming Zou}\\
 {\footnotesize\itshape Department of Mathematical Sciences, Tsinghua University, Beijing 100084, China.}}
\date{}
\begin{document}
\maketitle
\begin{abstract}
{\small
Let $n\ge2$ be an integer, let $1<p<n$, and set
$p^*=np/(n-p)$.
On the unit sphere $\mathbb S^n$ with its standard
Riemannian volume measure, we determine all powers
$\theta>0$ for which
\begin{equation*}
 \|u\|_{p^*}^{\theta}
 \le K_{n,p}^{\theta}\|\nabla u\|_p^{\theta}
      +B\|u\|_p^{\theta},
 \qquad u\in W^{1,p}(\mathbb S^n),
\end{equation*}
holds with some finite constant $B$, where $K_{n,p}$
is the sharp Euclidean Sobolev constant.
For $1<p<2$ and $n>p^2$, the largest admissible
power is $\min\{2,p^*,n/p\}$.
For $p>2$ and $3p-2<n<p^2$, it is
$2p(n-p)/(-3p^2+np+2n)$, completing the range
obtained by Aubin and Li~\cite{AL}.

For $\varepsilon\ge0$, let $B_\varepsilon(n,p)$
be the optimal lower-order coefficient in the
above inequality with $\theta=p$ and gradient
coefficient $K_{n,p}^p+\varepsilon$, with
$B_\varepsilon(n,p)=+\infty$ if no finite
coefficient is admissible.
For $p>2$ and $n\ge3p-2$, we determine its sharp
rate of divergence as $\varepsilon\downarrow0$,
including the logarithmic correction at $n=p^2$
and the logarithmic growth at $n=3p-2$.

We also obtain explicit bounds for the second
Sobolev constant $B_0(n,p)$.
Writing $V_n=|\mathbb S^n|$, we prove that
$B_0(n,p)=V_n^{-p/n}+O_n((2-p)^2)$ as $p\uparrow2$
for each fixed $n\ge3$, yielding its left derivative
at $p=2$.
In dimension four,
$\liminf_{p\downarrow2}B_0(4,p)\ge\frac43B_0(4,2)$,
which implies discontinuity at $p=2$.
For each fixed $n\ge5$, setting $P=(n+2)/3$,
we show that $B_0(n,p)\asymp_n(p-P)^{-P/2}$
as $p\downarrow P$.


}
\end{abstract}

\clearpage
\setcounter{tocdepth}{2}
\tableofcontents

\section{Introduction}\label{sec:introduction}

Let $n\ge2$ be an integer, let $1<p<n$, and put $p^*=np/(n-p)$.
The sharp Euclidean Sobolev inequality is
\begin{equation}\label{eq:euclidean-intro}
 \|v\|_{L^{p^*}(\R^n)}\le K_{n,p}\|\nabla v\|_{L^p(\R^n)},
 \qquad v\in\D^{1,p}(\R^n).
\end{equation}
Here $\D^{1,p}(\R^n)$ is the completion of $C_c^\infty(\R^n)$ with respect to the norm $\|\nabla\cdot\|_{L^p(\R^n)}$, and $K_{n,p}$ is the best constant in~\eqref{eq:euclidean-intro}. The classical results of Aubin~\cite{Aubin} and Talenti~\cite{Talenti} give
\begin{equation}\label{eq:euclidean-constant}
 K_{n,p}^p=\frac1n\left(\frac{p-1}{n-p}\right)^{p-1}
 \left[\frac{\Gamma(n+1)}
 {|\mathbb S^{n-1}|\,\Gamma(n/p)\Gamma(n+1-n/p)}\right]^{p/n}.
\end{equation}
The nonzero real-valued extremizers
of~\eqref{eq:euclidean-intro} are precisely
the functions of the form
\begin{equation}\label{eq:euclidean-extremals}
 v(x)=c\bigl(\lambda+|x-x_0|^{p/(p-1)}\bigr)^{-(n-p)/p},
 \qquad c\in\R\setminus\{0\},\quad
 \lambda>0,\quad x_0\in\R^n.
\end{equation}
See also Figalli--Zhang~\cite[Section 1.1]{FZ}.

On a smooth closed $n$-dimensional Riemannian manifold $(M,g)$, Aubin~\cite[Theorem 9, p.~593]{Aubin} proved that for every $\varepsilon>0$ there exists $C_\varepsilon<\infty$ such that
\begin{equation}\label{eq:compact-intro}
 \|u\|_{L^{p^*}(M)}^p
 \le (K_{n,p}^p+\varepsilon)\|\nabla u\|_{L^p(M)}^p
       +C_\varepsilon\|u\|_{L^p(M)}^p
 \qquad\bigl(u\in W^{1,p}(M)\bigr).
\end{equation}
Concentrating test functions in normal coordinates show that $K_{n,p}^p$ is the infimum of the admissible gradient coefficients; see also Hebey~\cite[Chapter 4 and Theorem 4.5, p.~95]{HebeyBook}. On the other hand, testing~\eqref{eq:compact-intro} with nonzero constant functions yields the lower bound $\operatorname{Vol}_g(M)^{-p/n}$ for the lower-order coefficient. The A-program asks whether the optimal gradient coefficient $K_{n,p}^p$ is admissible and, if so, seeks the smallest corresponding lower-order coefficient. The B-program reverses the order of optimization, fixing the lower-order coefficient at its volume lower bound and seeking the smallest admissible gradient coefficient, if one exists. See Hebey~\cite[Chapter 4]{HebeyBook} and Druet--Hebey~\cite{DH}.

Throughout this paper, $\mathbb S^n$ denotes the standard unit sphere with Riemannian volume measure $dV$, and $V_n=|\mathbb S^n|$. Unless otherwise specified, all integrals and $L^r$ norms on the sphere are taken with respect to $dV$. We write
\[
 \|u\|_r=\left(\int_{\mathbb S^n}|u|^r\,dV\right)^{1/r}.
\]

Aubin~\cite[Theorem 8, pp.~588--591]{Aubin} proved that the inequality
\[
 \|u\|_{p^*}^{\theta}\le K_{n,p}^{\theta}\|\nabla u\|_p^{\theta}+B\|u\|_p^{\theta},
 \qquad u\in W^{1,p}(\mathbb S^n),
\]
holds with some finite constant $B$ for $\theta=p$ if $1<p<2$, and for $\theta=p/(p-1)$ if $2\le p<n$. These results are also presented in Hebey~\cite[Propositions 4.3 and 4.6, pp.~99 and 114]{HebeyBook}.
Whenever this inequality holds, the gradient coefficient $K_{n,p}^{\theta}$ is optimal, as follows from the same concentration argument.

For $p=2$ and $n\ge3$, the sharp inequality on the sphere reads
\begin{equation}\label{eq:quadratic-intro}
 \|u\|_{2^*}^2\le K_{n,2}^2\|\nabla u\|_2^2
          +V_n^{-2/n}\|u\|_2^2,
 \qquad K_{n,2}^2=\frac{4V_n^{-2/n}}{n(n-2)},
 \qquad 2^*=\frac{2n}{n-2}.
\end{equation}
The gradient coefficient is optimal by concentration, and the lower-order coefficient is optimal by testing constant functions. This spherical inequality goes back to Aubin; see the account and the description of its extremizers in Hebey~\cite[Theorem 5.1, pp.~121--124]{HebeyBook}. Ilias~\cite[Theorem 3, p.~158]{Ilias1983} obtained the corresponding sharp comparison with the spherical inequality under a positive lower bound on the Ricci curvature. Beckner~\cite[Theorem 4 and (33)--(35), p.~232]{Beckner1993} proved a sharp family of Sobolev interpolation inequalities on the sphere for target exponents $2\le r\le 2^*$.

On general closed Riemannian manifolds, the sharp Sobolev inequality with outer power $\min\{p,2\}$ was established by Hebey--Vaugon~\cite{HV} for $p=2$, and independently by Aubin--Li~\cite[Theorem 1.1, p.~355]{AL} and Druet~\cite[p.~328]{Druet} for $p\ne2$.
If the sharp inequality holds for some power $\theta_0>0$, then it also holds for every $0<\theta\le\theta_0$, by the subadditivity of $t\mapsto t^{\theta/\theta_0}$.
For $p>2$, larger powers need not be admissible. If the scalar curvature of the manifold is positive somewhere, Druet~\cite[Theorem 5, p.~240]{Druet98} showed that the sharp inequality fails for every $\theta>2$ when $n>p^2$; see also~\cite[Section 6.2, p.~80]{DruetLecture2002}. Under the same curvature assumption, the sharp $p$-power inequality fails whenever $n>3p-2$; see~\cite[Remark 2, p.~114]{DruetNote}, \cite[Theorem 6.1, p.~76]{DruetLecture2002}, and \cite[Theorem 8.4, p.~76]{DH}.

In the intermediate range $p>2$, $3p-2<n<p^2$, Aubin--Li~\cite[Theorem 1.3, p.~356]{AL} proved the sharp inequality for every
\begin{equation}\label{eq:AL-range}
 0<\theta<\frac{2p(n-p)}{-3p^2+np+2n}.
\end{equation}
The upper endpoint in~\eqref{eq:AL-range} is not covered by their theorem. They also established the sharp inequality for every $0<\theta<p$ when $p>2$ and $n\le3p-2$.

The local Dirichlet problem on spherical caps is also relevant to the global sharp Sobolev inequality. Bandle--Brillard--Flucher~\cite{BBF1998} extended harmonic transplantation to spaces of constant curvature, with applications to the quadratic Dirichlet problem. On $\mathbb S^3$, Bandle--Peletier~\cite[Theorem 1, p.~84]{BandlePeletier1999} proved that the best constant in
\[
 \|u\|_6\le C\|\nabla u\|_2,\qquad u\in W_0^{1,2}(B_R),
\]
where $B_R$ is a geodesic ball of radius $R$, equals $K_{3,2}$ for $0<R\le\pi/2$, in which case there are no extremizers, whereas it is strictly larger than $K_{3,2}$ for $\pi/2<R<\pi$, and extremizers exist.
For general $1<p<n$, Bandle--Peletier--Stingelin~\cite[Theorem 2.1 and pp.~1389--1390]{BPS} proved that if
\[
 \frac{n-p}{p-1}<2,
\]
then on sufficiently small geodesic balls $B_R\subset\mathbb S^n$, the sharp Dirichlet Sobolev inequality
\[
 \|u\|_{p^*}\le K_{n,p}\|\nabla u\|_p,\qquad u\in W_0^{1,p}(B_R),
\]
holds. Combining this result with the local-to-global principle
of Druet~\cite[Lemma 6.2, p.~76]{DruetLecture2002}
and Aubin--Li~\cite[Theorem 7.1, pp.~379--380]{AL},
we obtain the sharp $p$-power inequality on the whole
sphere whenever $p>2$ and $n<3p-2$.

Let $B_0(n,p)$ denote the optimal lower-order coefficient in the sharp $p$-power inequality on $\mathbb S^n$, with $B_0(n,p)=+\infty$ if no finite coefficient exists. For $p>2$ and $n>3p-2$, Druet's result implies that $B_0(n,p)=+\infty$. In the critical case $n=3p-2$ with $p>2$, Frank--Jin--Wang--Xiong~\cite[Section 3, Step 5, pp.~16--17]{FJWX} considered a family $b_\lambda$ of concentrating Talenti functions cut off in a fixed geodesic ball and obtained, as $\lambda\downarrow0$,
\begin{equation*}
\begin{aligned}
 \|b_\lambda\|_{p^*}^p-K_{n,p}^p\|\nabla b_\lambda\|_p^p
 &=c\lambda^2\log(1/\lambda)+O(\lambda^2),\\
 \|b_\lambda\|_p^p&=O(\lambda^2),
\end{aligned}
\end{equation*}
where $c>0$. Dividing the first expression by $\|b_\lambda\|_p^p$ shows that $B_0(n,p)=+\infty$ also in this critical case. The preceding results therefore yield
\begin{equation}\label{eq:known-finiteness}
 B_0(n,p)<\infty\quad\Longleftrightarrow\quad n<3p-2
 \qquad(2<p<n).
\end{equation}
In particular, combining~\eqref{eq:known-finiteness} with Aubin--Li's theorem, we see that when $p>2$ and $n=3p-2$, the sharp inequality holds precisely for $0<\theta<p$.

The B-program also exhibits a qualitative difference between the cases $p=2$ and $p>2$. For $p>2$, no finite gradient coefficient is admissible when the lower-order coefficient is fixed at $V_n^{-p/n}$, as can be seen by perturbing a constant function~\cite[Proposition 4.1, pp.~83--84]{HebeyBook}; see also~\cite[Theorem 8.11, pp.~91--93]{DH}. Consequently, if $p>2$ and $B_0(n,p)<\infty$, then
$B_0(n,p)>V_n^{-p/n}$. Explicit lower bounds improving upon $V_n^{-p/n}$ are given in Appendix~\ref{app:bounds}.

In the quadratic case, Druet--Robert~\cite[Theorem 0.3, p.~747]{DruetRobert2001} studied the asymptotic profiles of sub-extremals on $\mathbb S^n$, $n\ge4$, equipped with a metric $g$ conformal to the standard metric but not isometric to it up to scaling. Their analysis concerns the limit in which the lower-order coefficient approaches the second Sobolev constant $B_0(g)$ from below, while the gradient coefficient remains optimal. In contrast, we fix the round metric and study the dependence of the optimal lower-order coefficient on $p$ and on the excess in the gradient coefficient.

A related problem is to improve the gradient
coefficient by imposing moment constraints.
Aubin~\cite{Aubin1979} showed that under the
first-moment constraint
\[
 \int_{\mathbb S^n}x|u|^{p^*}\,dV=0,
\]
the gradient coefficient in~\eqref{eq:compact-intro}
can be reduced to $2^{-p/n}K_{n,p}^p+\varepsilon$
for every $\varepsilon>0$;
see also Hebey~\cite[Theorem 6.3, p.~166]{HebeyBook}.
Hang--Wang~\cite[Theorem 1.1]{HangWang2022}
extended this result to higher-order moment constraints.
Under the first-moment constraint, Frank--Jin--Wang--Xiong~\cite[Theorem 1.1]{FJWX} showed that the excess $\varepsilon$ can be removed in the $p$-power inequality if and only if $1<p\le2$ or $2<p<n<3p-2$. They also established the corresponding sharp inequality with outer power $2$ and gradient coefficient $2^{-2/n}K_{n,p}^2$ for every $2\le p<n$.
In the quadratic case, Wei--Yan--Ye~\cite[Theorem 1.1]{WeiYanYe2026} further determined the optimal lower-order coefficient under the first-moment constraint for $n\ge4$.

In a previous work~\cite[Theorems 1.1, 1.2 and 1.6]{ZhouZouConstraints2026}, we established sharp improved Sobolev inequalities for $p=2$ under moment constraints of orders $m=1,2,3$. We also showed that the failure of the unconstrained sharp $p$-power inequality implies the failure of its sharp improved counterparts under moment constraints. When the unconstrained sharp $p$-power inequality holds, we obtained quantitative remainder estimates for the corresponding almost sharp improved inequalities. In the present paper, we focus on the unconstrained problem.

Our first main result determines the complete range of
admissible outer powers.

\begin{theorem}[Optimal outer powers]\label{thm:powers}
Let $n\ge2$ be an integer, let $1<p<n$, and put
$p^*=np/(n-p)$. Denote by $\mathcal T(n,p)$ the set of
all $\theta>0$ for which there exists a finite constant
$B=B(n,p,\theta)>0$ such that
\[
 \|u\|_{p^*}^\theta
 \le K_{n,p}^\theta\|\nabla u\|_p^\theta
       +B\|u\|_p^\theta
 \qquad\text{for all real-valued }
 u\in W^{1,p}(\mathbb S^n).
\]
Then
\[
 \mathcal T(n,p)=
 \begin{cases}
 (0,p],&1<p\le2,\ n\le p^2,\\
 (0,\min\{2,p^*,n/p\}],&1<p\le2,\ n>p^2,\\
 (0,2],&p>2,\ n\ge p^2,\\
 \displaystyle
 \left(0,\frac{2p(n-p)}{-3p^2+np+2n}\right],
       &p>2,\ 3p-2<n<p^2,\\
 (0,p),&p>2,\ n=3p-2,\\
 (0,p],&p>2,\ n<3p-2.
 \end{cases}
\]
\end{theorem}

In the range $1<p<2$ and $n>p^2$,
Theorem~\ref{thm:powers} shows that the optimal outer power
is strictly larger than $p$ and determines its exact value.
In the intermediate range $p>2$ and $3p-2<n<p^2$,
the theorem identifies the upper endpoint
in~\eqref{eq:AL-range} as the largest admissible power,
thereby completing the range obtained by Aubin--Li.

We next study the optimal lower-order coefficient as the
gradient coefficient approaches its sharp value.
For convenience, write
\begin{equation}\label{eq:global-notation}
 s=p^*,\qquad K=K_{n,p},\qquad \kappa=K^p,
\end{equation}
and, for $u\in W^{1,p}(\mathbb S^n)$, set
\begin{equation}\label{eq:global-quantities}
 N(u)=\|u\|_s^p,\qquad
 E(u)=\int_{\mathbb S^n}|\nabla u|^p\,dV,\qquad
 L(u)=\int_{\mathbb S^n}|u|^p\,dV.
\end{equation}
For $\varepsilon\ge0$, let $B_\varepsilon(n,p)$ denote
the smallest constant $B$ such that
\[
 N(u)\le(\kappa+\varepsilon)E(u)+BL(u)
 \qquad\text{for all real-valued }
 u\in W^{1,p}(\mathbb S^n),
\]
with the convention that $B_\varepsilon(n,p)=+\infty$
if no finite constant $B$ exists. Equivalently,
\begin{equation}\label{eq:Beps-def}
 B_\varepsilon(n,p)=
 \sup_{u\in W^{1,p}(\mathbb S^n)\setminus\{0\}}
 \frac{N(u)-(\kappa+\varepsilon)E(u)}{L(u)}.
\end{equation}

By~\eqref{eq:compact-intro} and the constant test,
$V_n^{-p/n}\le B_\varepsilon(n,p)<\infty$ for every $\varepsilon>0$.
Moreover, $B_\varepsilon(n,p)$ is nonincreasing in $\varepsilon$, and
\[
 B_\varepsilon(n,p)\uparrow B_0(n,p)
 \qquad\text{as }\varepsilon\downarrow0.
\]
In particular, $B_0(n,p)$ is the second Sobolev constant introduced above.
In the parameter range $p>2$ and $n\ge3p-2$, we know
from~\eqref{eq:known-finiteness} that $B_0(n,p)=+\infty$.
Our next result determines the sharp rate at which
$B_\varepsilon(n,p)$ diverges in this range.

\begin{theorem}[Sharp growth of the lower-order coefficient]
\label{thm:remainder}
Let $p>2$ and let $n\ge3p-2$ be an integer.
Then there exist constants $c,C>0$ and
$\varepsilon_0\in(0,1)$, depending only on $n,p$,
such that
\[
 c f_{n,p}(\varepsilon)
 \le B_\varepsilon(n,p)
 \le C f_{n,p}(\varepsilon)
 \qquad\text{for }0<\varepsilon<\varepsilon_0,
\]
where
\[
 f_{n,p}(\varepsilon)=
 \begin{cases}
 \varepsilon^{-(p-2)/2},
       &n>p^2\quad\text{(I)},\\
 \varepsilon^{-p(n-3p+2)/(2(n-p))},
       &3p-2<n<p^2\quad\text{(II)},\\
 \displaystyle
 \frac{\varepsilon^{-(p-2)/2}}{\log(1/\varepsilon)},
       &n=p^2\quad\text{(III)},\\
 [\log(1/\varepsilon)]^{p/2},
       &n=3p-2\quad\text{(IV)}.
 \end{cases}
\]
\end{theorem}

Throughout the paper, we write $f\asymp g$ if
$cg\le f\le Cg$ for some positive constants $c,C$
in the range under consideration.
Dependence on parameters other than those already fixed
is indicated by subscripts when necessary.

In Region I, the asserted growth rate can already be
deduced from the classical sharp Sobolev inequality
with outer power $2$ and the usual curvature expansion
for concentrating Talenti functions.
The logarithmic factors in Regions III and IV
cannot be recovered from the outer-power
classification and have different origins.
At $n=p^2$, the logarithm arises from the borderline
$L^p$-integrability of the Euclidean Talenti functions,
whereas at $n=3p-2$, it arises from the logarithmic
divergence of the second moment of their $p$-energy density.

We finally study the dependence of the second Sobolev constant on $p$.

Barbosa--Montenegro~\cite[Proposition 2.1, p.~969]{BM2010}
proved the continuity of the second Sobolev constant
on $[1,\min\{2,\sqrt n\})$ for general closed
Riemannian manifolds and raised the question of
continuity at $p=2$ from below.
On the round sphere, we give a quantitative answer
to this question for every $n\ge3$.

By Theorem~\ref{thm:powers} and~\eqref{eq:known-finiteness},
within the range $1<p<n$, one has
\[
 B_0(n,p)<\infty
 \quad\Longleftrightarrow\quad
 p\le2\ \text{or}\ p>\frac{n+2}{3}.
\]
The threshold $(n+2)/3$ coincides with $2$ precisely
when $n=4$.
Although $B_0(4,p)$ is finite for every $1<p<4$,
the following theorem shows that it is discontinuous
at $p=2$.

\pagebreak[3]
\begin{theorem}[Asymptotic behavior of the second Sobolev constant]
\label{thm:second}
The following assertions hold.
\begin{enumerate}
\renewcommand{\labelenumi}{(\alph{enumi})}
\item For each fixed $n\ge3$, one has
$B_0(n,2)=V_n^{-2/n}$ and
\[
 B_0(n,p)=V_n^{-p/n}+O_n((2-p)^2)
 \qquad\text{as }p\uparrow2.
\]
In particular, $p\mapsto B_0(n,p)$ is left-differentiable
at $p=2$, with
\[
 B'_{0,-}(n,2)=-\frac{\log V_n}{n}V_n^{-2/n}.
\]
\item In dimension four, one has
\[
 \liminf_{p\downarrow2}B_0(4,p)
 \ge\frac43V_4^{-1/2}>B_0(4,2).
\]
Consequently, $p\mapsto B_0(4,p)$ is discontinuous
at $p=2$.
\item For each fixed $n\ge5$, put $P=(n+2)/3$. Then
\[
 B_0(n,p)\asymp_n(p-P)^{-P/2}
 \qquad\text{as }p\downarrow P.
\]
Moreover,
\[
 \liminf_{p\downarrow P}(p-P)^{P/2}B_0(n,p)
 \ge\mathcal L_n>0,
\]
where $\mathcal L_n$ is the explicit constant
 given in~\eqref{eq:D-Ln}.
\end{enumerate}
\end{theorem}

In Appendix~\ref{app:bounds}, we also establish explicit upper and lower
bounds for $B_0(n,p)$ throughout its finite range. For $1<p\le2$, we
obtain an upper bound in terms of gamma functions and, for each fixed
$n\ge3$, an explicit bound on
$\limsup_{p\uparrow2}(V_n^{p/n}B_0(n,p)-1)/(2-p)^2$.
For $p>2$ and $n<3p-2$, we derive explicit lower bounds strictly above
$V_n^{-p/n}$, together with an explicit upper bound.

To prove Theorem~\ref{thm:remainder}, we first use spherical rearrangement
and a radial change of variables preserving the $p$-Dirichlet energy
to obtain local Sobolev inequalities on geodesic balls. A quantitative localization estimate preserves the gradient coefficient
of the local inequality and gives explicit control of the lower-order
coefficient in terms of the radius. This yields the upper bounds in
Regions II and IV. The upper bound in Region I follows from
averaging over the poles of stereographic projection. For the lower
bounds, we use cut-off Talenti functions, choosing the cutoff radius as
a function of the concentration scale in Regions II and IV to balance
the cutoff error against the curvature contribution.

For the critical upper bound at $n=p^2$, we use
a stability argument related to that of
Frank--Jin--Wang--Xiong~\cite[Section 2]{FJWX}.
We apply the gradient stability theorem of
Figalli--Zhang~\cite{FZ} to the Euclidean transplants
of cut-off extremizers; see also the earlier work
of Figalli--Neumayer~\cite{FigalliNeumayer2019}.
The key is to retain the negative contribution of
the Euclidean Sobolev deficit in the estimate for
the spherical deficit.
This yields sufficiently precise control of the error
in the approximation by Talenti functions.
Applying the reverse triangle inequality on a suitably
chosen shrinking Euclidean ball, we obtain a lower bound
of order $\lambda^p\log(1/\lambda)$ for the $L^p$ mass
of the spherical extremizers, where $\lambda$ is the
concentration scale of the approximating Talenti function.
Combining this lower bound with the deficit estimate
gives the critical upper bound.

For outer powers below $p$, we establish an equivalence
between admissibility and a corresponding growth bound
for $B_\varepsilon(n,p)$.
When $p>2$ and $n\ge3p-2$, the divergence of
$B_\varepsilon(n,p)$ excludes all powers $\theta\ge p$,
while the equivalence, together with
Theorem~\ref{thm:remainder}, determines the remaining
admissible powers.
For $1<p<2$ and $n>p^2$, the local inequalities
needed for the classification fall within the range
proved by Aubin--Li~\cite[Proposition 5.1]{AL}.
We derive these inequalities on the sphere
in Lemma~\ref{lem:C-local} and use height truncation
to extend them to the whole sphere for powers
exceeding $p$, beyond the scope of their
local-to-global theorem~\cite[Theorem 7.1]{AL}.
The remaining necessary bounds follow from concentrating
test functions, in some cases with an added constant.

For Theorem~\ref{thm:second}, the estimates as $p\uparrow2$
follow from averaging the stereographic energy identity
and expanding the resulting explicit upper bound in terms
of gamma functions. For $p>2$, we use the radial change of
variables together with estimates that are uniform in $p$
near the critical parameters. The quantitative localization
estimate yields the upper bound for $n\ge5$, while
transplanted Talenti functions give the corresponding
lower bound and the discontinuity result in dimension four.
In dimension four, it is crucial that the small-radius
expansion of the leading coefficient in the $L^p$-mass
asymptotics has a uniform relative error, even as
the associated beta integral diverges as $p\downarrow2$.

In Section~\ref{sec:preliminaries} we prove the common variational and geometric estimates. Section~\ref{sec:remainder} gives the growth of $B_\varepsilon$, and Section~\ref{sec:C-outer} proves the outer-power classification. Section~\ref{sec:D-second} treats the parameter limits of $B_0$. Appendix~\ref{app:bounds} contains the explicit bounds.

\section{Preliminaries and localization}\label{sec:preliminaries}

In addition to the notation introduced in Section~\ref{sec:introduction},
we set
\begin{equation}\label{eq:proof-notation}
 \nu=\frac ps=\frac{n-p}{n},\qquad
 q=\frac p{p-1},\qquad
 \beta=\frac{n-p}{p},\qquad
 a=\frac{n-p}{p-1},\qquad
 S=\kappa^{-1}.
\end{equation}
We write $\sigma=|\mathbb S^{n-1}|$.

We use the Talenti functions and the associated integrals
\begin{equation}\label{eq:talenti-data}
\begin{aligned}
 U(x)&=(1+|x|^q)^{-\beta},\qquad
 U_\lambda(x)=\lambda^{-\beta}U(x/\lambda),
 \quad \lambda>0,\\
 \mathsf A&=\int_{\R^n}U^s\,dx
 =\frac{\sigma}{q}\mathrm B(n/q,n/p),\\
 \mathsf E&=\int_{\R^n}|\nabla U|^p\,dx,
 \qquad \mathsf A^\nu=\kappa\mathsf E.
\end{aligned}
\end{equation}
By scaling, we have
\[
 \int_{\R^n}U_\lambda^s\,dx=\mathsf A,\qquad
 \int_{\R^n}|\nabla U_\lambda|^p\,dx=\mathsf E
 \quad\text{for every }\lambda>0.
\]
We also write $U(r)$ and $U_\lambda(r)$ for their radial
profiles.

Here $\mathrm B$ and $\Gamma$ denote the beta and gamma
functions, respectively. All integrals on $\R^n$ are taken
with respect to Lebesgue measure.

Unless otherwise specified, $B_R$ denotes the geodesic ball
of radius $R$ centered at a fixed pole $o\in\mathbb S^n$,
and $r=d(o,\cdot)$ denotes the geodesic distance from $o$.
In geodesic polar coordinates, the volume element is
$dV=(\sin r)^{n-1}\,dr\,d\omega$, where $d\omega$ is the
standard volume element on $\mathbb S^{n-1}$.

Throughout the proofs, $C$ and $c$ denote positive constants
that may change from line to line. Unless otherwise specified,
$n$ and $p$ are fixed. When $p$ or other parameters vary,
we specify the relevant uniformity of the estimates.

For $u\in W^{1,p}(\mathbb S^n)$, let $u^\star$ denote
the spherical decreasing rearrangement of $|u|$
about a prescribed pole.
Then $u^\star$ belongs to $W^{1,p}(\mathbb S^n)$,
is equimeasurable with $|u|$, and has a nonnegative,
nonincreasing radial profile.
In particular,
\[
 \|u^\star\|_z=\|u\|_z
 \quad\text{whenever }1\le z<\infty
 \text{ and }u\in L^z(\mathbb S^n).
\]
The spherical P\'olya--Szeg\H{o} inequality gives
\[
 E(u^\star)\le E(u);
\]
see Chiacchio--Ricciardi~\cite[Lemma 3.3, p.~184]{CR}.
If $0<R<\pi$ and $u\in W^{1,p}_0(B_R)$, we extend $u$
by zero to $\mathbb S^n$ and rearrange it about the
center of $B_R$. By equimeasurability and radial monotonicity,
$u^\star$ vanishes almost everywhere outside $B_R$. Since
$u^\star\in W^{1,p}(\mathbb S^n)$ and $\partial B_R$
is smooth, its restriction to $B_R$ belongs to
$W^{1,p}_0(B_R)$.

\subsection{Attainment and quantitative localization}

We first show that the infimum in the following
constrained minimization problem is attained whenever
it is strictly below the Euclidean Sobolev threshold $S$.
Using spherical rearrangement, we give
a direct concentration-compactness argument; cf.\
Lions~\cite[Lemma I.1, p.~158]{Lions}.

\begin{lemma}\label{lem:attainment}
For $\alpha>0$, let
\[
 \Lambda(\alpha)=\inf\left\{E(u)+\alpha L(u):
 u\in W^{1,p}(\mathbb S^n),\
 \int_{\mathbb S^n}|u|^s\,dV=1\right\}.
\]
If $\Lambda(\alpha)<S$, then for any prescribed pole
$o\in\mathbb S^n$, the infimum is attained by a
nonnegative function $u\in W^{1,p}(\mathbb S^n)$
that is radial and nonincreasing about $o$.
Moreover, $u$ satisfies
\begin{equation}\label{eq:euler}
\begin{aligned}
 &\int_{\mathbb S^n}|\nabla u|^{p-2}
       \langle\nabla u,\nabla\varphi\rangle\,dV
 +\alpha\int_{\mathbb S^n}u^{p-1}\varphi\,dV\\
 &\qquad=\Lambda(\alpha)
       \int_{\mathbb S^n}u^{s-1}\varphi\,dV
 \qquad\text{for all }\varphi\in W^{1,p}(\mathbb S^n).
\end{aligned}
\end{equation}
\end{lemma}

\begin{proof}
The Sobolev embedding gives $N(u)\le C(E(u)+L(u))$.
Since $\alpha>0$, it follows that $\Lambda(\alpha)>0$.
By spherical rearrangement, we may choose a nonnegative
minimizing sequence $(u_j)$, radial and nonincreasing
about $o$, with
\[
 \int_{\mathbb S^n}u_j^s\,dV=1,\qquad
 E(u_j)+\alpha L(u_j)\longrightarrow\Lambda(\alpha).
\]
Since $\alpha>0$, the sequence is bounded in
$W^{1,p}(\mathbb S^n)$.
By reflexivity and the compact embedding
$W^{1,p}(\mathbb S^n)\hookrightarrow L^p(\mathbb S^n)$,
we may pass to a subsequence such that
\[
 u_j\rightharpoonup u\quad\text{in }W^{1,p}(\mathbb S^n),
 \qquad
 u_j\longrightarrow u\quad\text{in }L^p(\mathbb S^n)
 \text{ and a.e.}
\]
Passing to a further subsequence if necessary,
we also have
\[
 |\nabla u_j|^p\,dV\stackrel{*}{\rightharpoonup}\mu
 \quad\text{as finite measures on }\mathbb S^n.
\]
The limit $u$ is also nonnegative, radial and nonincreasing
about $o$. By Fatou's lemma,
\[
 M:=\int_{\mathbb S^n}u^s\,dV\in[0,1].
\]

For each fixed $r\in(0,\pi)$, radial monotonicity
and the normalization give
\[
 u_j(x)\le |B_r|^{-1/s}
 \qquad\text{for a.e. }x\in\mathbb S^n\setminus B_r.
\]
The same bound holds for $u$.
Hence the strong $L^p$ convergence implies that
$u_j\to u$ strongly in
$L^s(\mathbb S^n\setminus B_r)$.

For sufficiently small $\delta>0$, choose a Lipschitz
cutoff $0\le\chi_\delta\le1$ equal to one on $B_\delta$,
compactly supported in $B_{2\delta}$, and satisfying
$|\nabla\chi_\delta|\le C/\delta$. Then
\[
 \liminf_{j\to\infty}
 \int_{\mathbb S^n}|\chi_\delta u_j|^s\,dV
 \ge 1-\int_{\mathbb S^n\setminus B_\delta}u^s\,dV
 \ge 1-M.
\]
In normal coordinates centered at $o$, the metric
coefficients satisfy $g_{ij}=\delta_{ij}+O(\delta^2)$
on $B_{2\delta}$.
In these coordinates, extend $\chi_\delta u_j$
by zero to $\R^n$.
The sharp Euclidean Sobolev inequality then gives
\[
 (1-M)^\nu
 \le \kappa(1+C\delta^2)
 \liminf_{j\to\infty}
 \int_{\mathbb S^n}|\nabla(\chi_\delta u_j)|^p\,dV.
\]
For any $\xi>0$,
\[
\begin{aligned}
 \int_{\mathbb S^n}|\nabla(\chi_\delta u_j)|^p\,dV
 \le{}&(1+\xi)\int_{B_{2\delta}}|\nabla u_j|^p\,dV\\
 &+C_{p,\xi}\delta^{-p}
 \int_{B_{2\delta}}u_j^p\,dV.
\end{aligned}
\]
Letting $j\to\infty$ with $\delta$ and $\xi$ fixed,
and using strong $L^p$ convergence and the Portmanteau
theorem for the closed ball $\overline B_{2\delta}$,
we obtain
\[
\begin{aligned}
 (1-M)^\nu
 \le\kappa(1+C\delta^2)\biggl[
 &(1+\xi)\mu(\overline B_{2\delta})\\
 &+C_{p,\xi}\delta^{-p}
 \int_{B_{2\delta}}u^p\,dV
 \biggr].
\end{aligned}
\]
By H\"older's inequality,
\[
 \delta^{-p}\int_{B_{2\delta}}u^p\,dV
 \le C\left(\int_{B_{2\delta}}u^s\,dV\right)^{p/s}
 \longrightarrow0
 \qquad\text{as }\delta\downarrow0.
\]
Since $\mu(\overline B_{2\delta})\to\mu(\{o\})$,
letting $\delta\downarrow0$ and then $\xi\downarrow0$
gives
\begin{equation}\label{eq:atomic-cost}
 \mu(\{o\})\ge S(1-M)^\nu.
\end{equation}

By homogeneity, the definition of $\Lambda(\alpha)$ gives
$E(v)+\alpha L(v)\ge\Lambda(\alpha)N(v)$ for every
$v\in W^{1,p}(\mathbb S^n)$.
By local weak lower semicontinuity,
$\mu\ge|\nabla u|^p\,dV$ in the sense of measures.
Since $|\nabla u|^p\,dV$ has no atom at $o$,
combining~\eqref{eq:atomic-cost} with strong $L^p$
convergence and the preceding homogeneity inequality
yields
\begin{equation*}
 \Lambda(\alpha)=\mu(\mathbb S^n)+\alpha L(u)
 \geq E(u)+\alpha L(u)+S(1-M)^\nu
 \geq\Lambda(\alpha)M^\nu+S(1-M)^\nu\,.
\end{equation*}

Since $0<\nu<1$, we have
$M^\nu\ge M$ and $(1-M)^\nu\ge1-M$. Hence
\[
 \Lambda(\alpha)\ge
 \Lambda(\alpha)+(S-\Lambda(\alpha))(1-M).
\]
As $\Lambda(\alpha)<S$, we conclude that $M=1$.
Weak lower semicontinuity then shows that $u$
attains the infimum.

Equivalently, $u$ minimizes the quotient
\[
 \frac{E(v)+\alpha L(v)}{N(v)}
 \qquad\text{over }
 v\in W^{1,p}(\mathbb S^n)\setminus\{0\}.
\]
Differentiating this quotient at $u$ in any direction
$\varphi\in W^{1,p}(\mathbb S^n)$ and using
$N(u)=1$ and $E(u)+\alpha L(u)=\Lambda(\alpha)$,
we obtain~\eqref{eq:euler}.
\end{proof}

\begin{lemma}[Exact duality]\label{lem:duality}
For $A>\kappa$, set $B=B_{A-\kappa}(n,p)$. Then
\[
 \Lambda(B/A)=A^{-1}<S.
\]
Moreover, for any prescribed pole $o\in\mathbb S^n$,
there exists a nonnegative function
$u\in W^{1,p}(\mathbb S^n)$ that is radial and
nonincreasing about $o$, with
$\int_{\mathbb S^n}u^s\,dV=1$, such that
\begin{equation}\label{eq:extremizer}
 AE(u)+BL(u)=1,\qquad
 -\Delta_pu+\frac BAu^{p-1}=A^{-1}u^{s-1},
\end{equation}
where the equation is understood in the weak sense
on $\mathbb S^n$ and
$\Delta_pu=\operatorname{div}
(|\nabla u|^{p-2}\nabla u)$.
\end{lemma}

\begin{proof}
Since $A>\kappa$, we have $0<B<\infty$.
By the definition of $B$,
\[
 N(v)\le AE(v)+BL(v)
 \qquad\text{for all }v\in W^{1,p}(\mathbb S^n).
\]
In particular, $\Lambda(B/A)\ge A^{-1}$.

Suppose that $\ell:=\Lambda(B/A)>A^{-1}$.
By homogeneity, for every $v\in W^{1,p}(\mathbb S^n)$,
\[
 N(v)\le\ell^{-1}E(v)+\frac{B}{A\ell}L(v)
 \le AE(v)+\frac{B}{A\ell}L(v).
\]
Since $B/(A\ell)<B$, this contradicts the optimality
of $B$. Hence $\Lambda(B/A)=A^{-1}<S$.

Applying Lemma~\ref{lem:attainment} with $\alpha=B/A$,
we obtain a nonnegative minimizer $u$ that is radial
and nonincreasing about the prescribed pole $o$.
Its normalization, the identity $\Lambda(B/A)=A^{-1}$,
and the Euler--Lagrange equation~\eqref{eq:euler}
yield~\eqref{eq:extremizer}.
\end{proof}

We next establish a quantitative version of the
local-to-global principle discussed in the introduction,
with explicit dependence on the radius of the geodesic ball.

\begin{proposition}[Quantitative localization]
\label{prop:localization}
Let $0<R\le1$ and $A\ge\kappa$. Suppose that
\[
 N(f)\le AE(f)
 \qquad\text{for every }f\in W^{1,p}_0(B_R).
\]
Set
\[
 m_R=|B_{R/4}|^{-p/n},\qquad
 d_R=\frac4R,\qquad H_p=[2(p-1)]^{p-1},
\]
and
\[
 T_R=2(Sm_R+H_pd_R^p),\qquad
 Q_R=(2^{p-1}-1)T_R+2^{p-1}d_R^p.
\]
Then
\begin{equation}\label{eq:localization}
 N(u)\le AE(u)+(m_R+AQ_R)L(u)
 \qquad\text{for every }u\in W^{1,p}(\mathbb S^n).
\end{equation}
Moreover, for every compact interval $I\Subset(1,n)$,
there exists a constant $C(n,I)>0$ such that
\[
 m_R+AQ_R\le C(n,I)AR^{-p}
\]
uniformly for $p\in I$, $0<R\le1$, and
$A\ge\kappa(n,p)$.

In particular, if the local inequality holds with
$A=\kappa$, then
\begin{equation}\label{eq:explicit-upper}
 B_0(n,p)\le
 (2^p-1)m_R+\kappa C_p4^pR^{-p},
 \qquad
 C_p=2(2^{p-1}-1)H_p+2^{p-1}.
\end{equation}
\end{proposition}

\begin{proof}
We first establish cutoff estimates for normalized
radial weak solutions, without using the assumed
local Sobolev inequality. Let $u\in W^{1,p}(\mathbb S^n)$
be nonnegative, radial and nonincreasing about the center
of $B_R$, with $\int_{\mathbb S^n}u^s\,dV=1$, and suppose
that
\[
 -\Delta_pu+\alpha u^{p-1}=\lambda u^{s-1}
\]
weakly on $\mathbb S^n$, where $\alpha\ge0$ and
$0\le\lambda\le S$.

By radial monotonicity and normalization,
\[
 u(r)^{s-p}\le |B_r|^{-p/n}\le m_R
 \qquad\text{for a.e. }r\in[R/4,\pi).
\]
Choose a radial Lipschitz cutoff $0\le\zeta\le1$
that vanishes on $B_{R/4}$, equals one outside
$B_{R/2}$, and satisfies $|\nabla\zeta|\le d_R$.

Testing the equation with
$\zeta^pu\in W^{1,p}(\mathbb S^n)$, discarding
the nonnegative zeroth-order term, and using
$\lambda\le S$ and $u^{s-p}\le m_R$ on the support
of $\zeta$, we obtain
\[
\begin{aligned}
 \int_{\mathbb S^n}\zeta^p|\nabla u|^p\,dV
 \le{}&Sm_RL(u)\\
 &+p\int_{\mathbb S^n}
 \zeta^{p-1}u|\nabla u|^{p-1}
 |\nabla\zeta|\,dV.
\end{aligned}
\]
For $X,Y\ge0$, Young's inequality gives
\[
 pX^{p-1}Y\le\frac12X^p+H_pY^p.
\]
Applying this with $X=\zeta|\nabla u|$ and
$Y=u|\nabla\zeta|$, and absorbing the first term,
we obtain
\begin{equation}\label{eq:annular-energy}
 \int_{\mathbb S^n\setminus B_{R/2}}
 |\nabla u|^p\,dV
 \le 2(Sm_R+H_pd_R^p)L(u)
 =T_RL(u).
\end{equation}

Let $0\le\eta\le1$ be any radial nonincreasing Lipschitz
cutoff that equals one on $B_{R/2}$, vanishes outside
$B_{3R/4}$, and satisfies $|\nabla\eta|\le d_R$.
Since $\eta u$ has compact support in $B_R$, it belongs
to $W^{1,p}_0(B_R)$.

Since $\eta=1$ on $B_{R/2}$, the energy difference
$E(\eta u)-E(u)$ comes entirely from
$\mathbb S^n\setminus B_{R/2}$.
Applying $|X+Y|^p\le2^{p-1}(|X|^p+|Y|^p)$ to
$\nabla(\eta u)=\eta\nabla u+u\nabla\eta$ there
and using~\eqref{eq:annular-energy}, we obtain
\[
\begin{aligned}
 E(\eta u)-E(u)
 &\le (2^{p-1}-1)
 \int_{\mathbb S^n\setminus B_{R/2}}
 |\nabla u|^p\,dV\\
 &\quad+2^{p-1}d_R^pL(u)
 \le Q_RL(u).
\end{aligned}
\]
The corresponding lower bound follows by discarding
the energy of $\eta u$ outside $B_{R/2}$ and
using~\eqref{eq:annular-energy}.
Moreover, the pointwise bound $u^{s-p}\le m_R$
outside $B_{R/2}$ gives
\[
 \int_{\mathbb S^n\setminus B_{R/2}}u^s\,dV
 \le m_RL(u).
\]
Consequently,
\begin{equation}\label{eq:cutoff-two-sided}
\begin{aligned}
 -T_RL(u)&\le E(\eta u)-E(u)\le Q_RL(u),\\
 0&\le1-\int_{\mathbb S^n}(\eta u)^s\,dV
 \le m_RL(u).
\end{aligned}
\end{equation}
The constants in these estimates are independent
of $\alpha$ and $\lambda$ in the specified ranges.

We now prove~\eqref{eq:localization}.
Set $\alpha=Q_R+m_R/A$, so that
$A\alpha=m_R+AQ_R$.
By homogeneity and the definition of $\Lambda(\alpha)$,
it suffices to prove that $\Lambda(\alpha)\ge A^{-1}$.

Suppose instead that $\Lambda(\alpha)<A^{-1}\le S$.
By Lemma~\ref{lem:attainment}, there exists a normalized
nonnegative minimizer $u$ that is radial and nonincreasing
about the center of $B_R$. It satisfies
\[
 -\Delta_pu+\alpha u^{p-1}
 =\Lambda(\alpha)u^{s-1}.
\]
Since $\alpha>0$ and $0<\Lambda(\alpha)<A^{-1}\le S$,
the cutoff estimates~\eqref{eq:cutoff-two-sided}
apply to $u$.

Since $0<\nu<1$ and
$\int_{\mathbb S^n}(\eta u)^s\,dV\le1$, we have
\[
 N(\eta u)\ge
 \int_{\mathbb S^n}(\eta u)^s\,dV
 \ge1-m_RL(u).
\]
Applying the local inequality to $\eta u$, we obtain
\[
\begin{aligned}
 1-m_RL(u)
 &\le N(\eta u)
 \le AE(\eta u)\\
 &\le A[E(u)+Q_RL(u)]\\
 &=A\Lambda(\alpha)-m_RL(u)\\
 &<1-m_RL(u),
\end{aligned}
\]
a contradiction. Thus $\Lambda(\alpha)\ge A^{-1}$,
and~\eqref{eq:localization} follows.

It remains to verify the uniform and explicit bounds.
Fix a compact interval $I\Subset(1,n)$.
Since $|B_{R/4}|\ge c_nR^n$ for $0<R\le1$,
and $\kappa$, $S$ and $H_p$ are continuous and
positive for $p\in I$, we have
\[
 m_R+T_R+Q_R\le C(n,I)R^{-p}.
\]
Since $\kappa(n,p)$ is bounded away from zero
uniformly for $p\in I$, the condition
$A\ge\kappa(n,p)$ implies
\[
 m_R+AQ_R\le C(n,I)AR^{-p}.
\]
Finally, substituting $S=\kappa^{-1}$ into
$m_R+\kappa Q_R$ gives
\[
 m_R+\kappa Q_R
 =(2^p-1)m_R+\kappa C_p4^pR^{-p},
\]
which proves~\eqref{eq:explicit-upper}.
\end{proof}

\subsection{Powers and remainders}

\begin{lemma}\label{lem:power}
For fixed $n,p$, let $0<\theta<p$ and set
\[
 \rho=\frac{\theta}{p},\qquad
 b=\frac{1-\rho}{\rho}=\frac p\theta-1.
\]
Then there exists a finite constant $C_0$ such that
\begin{equation}\label{eq:power}
 N(u)^\rho\le\kappa^\rho E(u)^\rho+C_0L(u)^\rho
 \qquad\text{for all }u\in W^{1,p}(\mathbb S^n)
\end{equation}
if and only if
\begin{equation}\label{eq:growth}
 B_\varepsilon(n,p)=O(\varepsilon^{-b})
 \qquad\text{as }\varepsilon\downarrow0.
\end{equation}
Moreover, whenever~\eqref{eq:power} holds with
a constant $C_0$, we have
\begin{equation}\label{eq:envelope}
 B_\varepsilon(n,p)\le
 C_0^{1/\rho}
 \left[1-(1+\varepsilon/\kappa)^{-1/b}\right]^{-b}
 \qquad\text{for every }\varepsilon>0.
\end{equation}
The constant $C_0$ and the implicit constant
in~\eqref{eq:growth} may depend on $n,p,\theta$.
\end{lemma}

\begin{proof}
Suppose that~\eqref{eq:power} holds. Testing it with
constant functions gives $C_0\ge V_n^{-\theta/n}>0$.
For $x,y\ge0$ and $0<\tau<1$, applying the concavity
of $t\mapsto t^\rho$ with weights $\tau$ and $1-\tau$
to $x/\tau^{1/\rho}$ and $y/(1-\tau)^{1/\rho}$ gives
\[
 (x^\rho+y^\rho)^{1/\rho}
 \le\tau^{-b}x+(1-\tau)^{-b}y.
\]
For $\varepsilon>0$, choose
$\tau=(1+\varepsilon/\kappa)^{-1/b}$.
Applying this inequality with $x=\kappa E(u)$ and
$y=C_0^{1/\rho}L(u)$, we obtain from~\eqref{eq:power}
\[
 N(u)\le(\kappa+\varepsilon)E(u)
       +C_0^{1/\rho}(1-\tau)^{-b}L(u).
\]
By the definition of $B_\varepsilon(n,p)$, this proves
\eqref{eq:envelope}. Since
$1-\tau\sim\varepsilon/(b\kappa)$ as
$\varepsilon\downarrow0$, we obtain~\eqref{eq:growth}.

Conversely, suppose that~\eqref{eq:growth} holds.
By the monotonicity of $B_\varepsilon(n,p)$ in
$\varepsilon$, this estimate extends, with a
possibly larger constant $C>0$, to
\[
 B_\varepsilon(n,p)\le C\varepsilon^{-b}
 \qquad(0<\varepsilon\le1).
\]
For any nonzero $u$, set $t=E(u)/L(u)$ and choose
$\varepsilon=(1+t)^{-\rho}$. Since
$\rho b=1-\rho$, the definition of $B_\varepsilon$
gives
\[
\begin{aligned}
 \frac{N(u)}{L(u)}
 &\le\kappa t+\varepsilon t+C\varepsilon^{-b}\\
 &\le\kappa t+(1+C)(1+t)^{1-\rho}.
\end{aligned}
\]
Replacing $\kappa t$ by $\kappa(1+t)$ in the
preceding upper bound, we apply the tangent-line
inequality for $z\mapsto z^\rho$ at
$z=\kappa(1+t)$.
Using also the subadditivity of this function,
we obtain
\[
\begin{aligned}
 \left(\frac{N(u)}{L(u)}\right)^\rho
 &\le\kappa^\rho(1+t)^\rho
       +\rho\kappa^{\rho-1}(1+C)\\
 &\le\kappa^\rho t^\rho+\kappa^\rho
       +\rho\kappa^{\rho-1}(1+C).
\end{aligned}
\]
Multiplying by $L(u)^\rho$ proves~\eqref{eq:power}.
The case $u=0$ is immediate.
\end{proof}

By definition, $p\in\mathcal T(n,p)$ if and only if
$B_0(n,p)<\infty$.
Since admissibility at a given power implies
admissibility at every smaller positive power,
the finiteness of $B_0(n,p)$ is also necessary
for the admissibility of any $\theta>p$.
Powers above $p$ are treated separately in
Section~\ref{sec:C-outer}.

\subsection{Radial transplantation}\label{sec:transplant}

We use a radial version of harmonic transplantation
adapted to the $p$-Dirichlet energy; see
Bandle--Brillard--Flucher~\cite{BBF1998} for the
quadratic case on spaces of constant curvature,
and Sano--Takahashi~\cite{ST2022} for related
$p$-harmonic transformations in the Euclidean setting.

For $0<R<\pi$ and $0<r<R$, define
\begin{align}\label{eq:transplant-def}
 G_R(r)&=a\int_r^R(\sin\tau)^{-a-1}\,d\tau,
 \qquad t(r)=G_R(r)^{-1/a},\nonumber\\
 H_R(r)&=(\sin r)^{-a}-G_R(r),
 \qquad F_R(r)=(\sin r)^aG_R(r).
\end{align}
Since $G_R(r)\sim r^{-a}$ as $r\downarrow0$ and
$G_R(r)\to0$ as $r\uparrow R$, the map
$t:(0,R)\to(0,\infty)$ is an increasing
diffeomorphism. It satisfies
\[
 \lim_{r\downarrow0}\frac{t(r)}r=1,\qquad
 t'(r)=\left(\frac{t(r)}{\sin r}\right)^{a+1}.
\]
We denote its inverse by $r=r(t)$.

First consider radial functions
$w\in C_c^\infty(B_R\setminus\{o\})$ and
$v\in C_c^\infty(\R^n\setminus\{0\})$
related by $w(r)=v(t(r))$.
Since $(a+1)(p-1)=n-1$, we have
\[
 (t')^{p-1}(\sin r)^{n-1}=t^{n-1}.
\]
Identifying $w$ with its zero extension to $\mathbb S^n$,
we obtain by a change of variables, for every $z>0$,
\begin{align}\label{eq:transplant}
 E(w)&=\int_{\R^n}|\nabla v|^p\,dx,\nonumber\\
 \int_{B_R}|w|^z\,dV
 &=\sigma\int_0^\infty
 |v(t)|^zW_R(t)t^{n-1}\,dt,
\end{align}
where the weight $W_R$ is given by
\[
 W_R(t)=\left(\frac{\sin r(t)}t\right)^{n+a}
       =F_R(r(t))^{(n+a)/a}.
\]

The correspondence between these two classes of
smooth functions is bijective.
Both classes are dense in the respective radial
Sobolev spaces: rotational averaging gives radial
smooth approximations, and since $p<n$, these
approximations can be cut off near the center with
arbitrarily small error in the $L^p$ norm of the gradient.
By Poincar\'e's inequality on $B_R$ and the definition
of $\D^{1,p}(\R^n)$, the energy identity extends this
correspondence to an isometric isomorphism between
the radial subspaces of $W^{1,p}_0(B_R)$ and
$\D^{1,p}(\R^n)$, both endowed with the gradient norm.

By the Sobolev embeddings, convergence in the
gradient norms implies convergence in $L^s$ on
both domains. Passing to subsequences and using
the fact that $t(r)$ is a diffeomorphism on $(0,R)$,
we recover the relation $w(r)=v(t(r))$ almost everywhere
for the limiting functions.
Consequently, both identities in~\eqref{eq:transplant}
remain valid, with both sides of the second identity
allowed to take the value $+\infty$.

For a general $f\in W^{1,p}_0(B_R)$, we first extend
$f$ by zero to $\mathbb S^n$, rearrange it about $o$,
and then transplant $f^\star|_{B_R}$ to $\R^n$.

\begin{lemma}\label{lem:weight}
Fix $a\ge2$. There exist constants $C,c>0$,
depending only on $a$, such that for every $0<R\le1$,
\[
 \sup_{0<r<R}(F_R(r)-1)_+\le\omega_R,
\]
where
\[
 \omega_R=
 \begin{cases}
 CR^{2a/(a-2)},&a>2,\\
 CR^2e^{-c/R^2},&a=2.
 \end{cases}
\]
Consequently, every $f\in W^{1,p}_0(B_R)$ satisfies
\begin{equation}\label{eq:local-A}
 N(f)\le
 \kappa(1+\omega_R)^{p(n-1)/n}E(f).
\end{equation}
\end{lemma}

\begin{proof}
For every $a>0$, differentiating
\eqref{eq:transplant-def} gives
\[
 H_R'(r)=a(1-\cos r)(\sin r)^{-a-1}.
\]
Integrating from $r$ to $R$ and using
$G_R(r)\to0$ as $r\uparrow R$, we obtain
\[
 H_R(r)=(\sin R)^{-a}
 -a\int_r^R(1-\cos\tau)(\sin\tau)^{-a-1}\,d\tau.
\]
We now assume $a\ge2$. Using
\[
 \frac56r\le\sin r\le r,\qquad
 1-\cos r\le\frac12r^2
 \qquad(0<r\le1),
\]
together with $F_R-1=-(\sin r)^aH_R$, we obtain
\[
 F_R(r)-1
 \le C_ar^a\int_r^R\tau^{1-a}\,d\tau
       -c_a(r/R)^a.
\]

If $a>2$, this yields
\[
 F_R(r)-1\le C_ar^2-c_aR^{-a}r^a.
\]
The representation of $H_R$ also gives
$|H_R(r)|\le C_a(R^{-a}+r^{2-a})$ in this case.
Hence, for fixed $R$, we have
$F_R(r)=1+O_{a,R}(r^2)$ as $r\downarrow0$.
Setting $r=R^{a/(a-2)}y$, we deduce that
\[
 F_R(r)-1
 \le R^{2a/(a-2)}(C_ay^2-c_ay^a).
\]
Since $C_ay^2-c_ay^a$ is bounded above for $y>0$,
we obtain the asserted bound on
$\sup_{0<r<R}(F_R(r)\allowbreak-1)_+$ in the case $a>2$.

If $a=2$, set $x=r/R$. Then
\[
 F_R(r)-1
 \le x^2[C_*R^2\log(1/x)-c_*].
\]
Here $C_*,c_*>0$ are independent of $R$.
The maximum of the right-hand side over $0<x<1$ is
\[
 \frac{C_*R^2}{2e}
 \exp\left(-\frac{2c_*}{C_*R^2}\right),
\]
which proves the estimate in this case.

Finally, let $v$ be the Euclidean transplant
of $f^\star$. Since
$W_R(t)\le(1+\omega_R)^{(n+a)/a}$,
combining~\eqref{eq:transplant} with the sharp
Euclidean Sobolev inequality and the spherical
P\'olya--Szeg\H{o} inequality yields
\[
\begin{aligned}
 N(f)=N(f^\star)
 &\le(1+\omega_R)^{\nu(n+a)/a}
      \|v\|_{L^s(\R^n)}^p\\
 &\le\kappa(1+\omega_R)^{p(n-1)/n}E(f^\star)\\
 &\le\kappa(1+\omega_R)^{p(n-1)/n}E(f),
\end{aligned}
\]
where we used $\nu(n+a)/a=p(n-1)/n$.
This proves~\eqref{eq:local-A}.
\end{proof}

\subsection{Stereographic projection and averaging}\label{sec:stereographic}

We use the stereographic averaging argument of
Frank--Jin--Wang--Xiong~\cite[Section 2, (5)--(8)]{FJWX}.
For later use, we compute explicitly the moments of the
vector field arising in the stereographic transformation.

\begin{lemma}\label{lem:stereographic}
Let $1<p<n$ and let $u\in C^\infty(\mathbb S^n)$
be real-valued. For $o\in\mathbb S^n$, denote by
$\Phi_o:\R^n\to\mathbb S^n\setminus\{-o\}$
the inverse stereographic projection from $-o$,
with $\Phi_o(0)=o$. Set
\[
 \rho_o(x)=\frac{2}{1+|x|^2}
           =1+o\cdot\Phi_o(x),
 \qquad
 v_o(x)=\rho_o(x)^\beta u(\Phi_o(x)),
\]
and define
\[
 Z_o(y)=\nabla_{\mathbb S^n}\log(1+o\cdot y),
 \qquad y\in\mathbb S^n\setminus\{-o\}.
\]
For convenience, set $Z_o(-o)=0$.
Then $v_o\in\D^{1,p}(\R^n)$ and
\begin{equation}\label{eq:D-stereo}
\begin{aligned}
 \int_{\R^n}|v_o|^s\,dx
 &=\int_{\mathbb S^n}|u|^s\,dV,\\
 \int_{\R^n}|\nabla v_o|^p\,dx
 &=\int_{\mathbb S^n}
 |\nabla_{\mathbb S^n}u+\beta uZ_o|^p\,dV.
\end{aligned}
\end{equation}
Let $\mathbb E_o$ denote expectation with respect to
the normalized volume measure $dV(o)/V_n$
on $\mathbb S^n$.
Then, for every fixed $y\in\mathbb S^n$,
\begin{equation}\label{eq:stereo-moments}
 \mathbb E_o Z_o(y)=0,\qquad
 \mathbb E_o|Z_o(y)|^p=M_{n,p}:=
 \frac{\Gamma((n+p)/2)\Gamma((n-p)/2)}
      {\Gamma(n/2)^2}.
\end{equation}
Moreover, if $n>2$, then
\[
 \mathbb E_o|Z_o(y)|^2=\frac{n}{n-2}.
\]
\end{lemma}

\begin{proof}
Fix $y\in\mathbb S^n$. For $o\notin\{y,-y\}$,
write $o=(\cos r)y+(\sin r)e$, where
$r=d(o,y)\in(0,\pi)$ and $e$ is a unit vector
in $T_y\mathbb S^n$.
\[
 Z_o(y)=\frac{o-(o\cdot y)y}{1+o\cdot y}
       =\tan(r/2)e.
\]
Conditional on $r$, the direction $e$ is uniformly
distributed on the unit sphere in $T_y\mathbb S^n$.
In particular, $Z_o(y)$ and $-Z_o(y)$ have the same
distribution.
For $0<m<n$, the substitution $t=\sin^2(r/2)$ gives
\[
\begin{aligned}
 \mathbb E_o|Z_o(y)|^m
 &=\frac{\int_0^\pi\tan^m(r/2)\sin^{n-1}r\,dr}
         {\int_0^\pi\sin^{n-1}r\,dr}\\
 &=\frac{\mathrm B((n+m)/2,(n-m)/2)}
         {\mathrm B(n/2,n/2)}.
\end{aligned}
\]
Taking $m=1$ shows that $Z_o(y)$ is integrable
with respect to $o$. By symmetry,
$\mathbb E_o Z_o(y)=0$.
Taking $m=p$ and, if $n>2$, $m=2$ gives the
remaining moment identities.
Moreover, since $|Z_o(y)|=|Z_y(o)|$, we have
\[
 \int_{\mathbb S^n}|Z_o(y)|^p\,dV(y)
 =V_nM_{n,p}<\infty
\]
for every fixed $o\in\mathbb S^n$.

The pullback of the round metric under $\Phi_o$
is $\rho_o^2g_{\R^n}$, with volume element
$\rho_o^n\,dx$.
By the chain rule, writing $y=\Phi_o(x)$, we have
\[
 |\nabla_xv_o(x)|
 =\rho_o(x)^{\beta+1}
 |\nabla_{\mathbb S^n}u(y)+\beta u(y)Z_o(y)|.
\]
Since $\beta s=n$ and $p(\beta+1)=n$, changing
variables gives both identities in~\eqref{eq:D-stereo}.
Since $u$ is smooth and $Z_o$ is $p$-integrable
on $\mathbb S^n$ for every fixed $o$, both integrals
on the right-hand side of~\eqref{eq:D-stereo}
are finite. Thus $v_o\in L^s(\R^n)$ and
$\nabla v_o\in L^p(\R^n)$.

For $T>0$, let $B_T$ denote the Euclidean ball
of radius $T$ centered at the origin.
Choose $\chi_T\in C_c^\infty(\R^n)$ with
$0\le\chi_T\le1$, equal to one on $B_T$,
supported in $B_{2T}$, and satisfying
$|\nabla\chi_T|\le C/T$.
By H\"older's inequality and $1/p-1/s=1/n$,
\[
\begin{aligned}
 \|\nabla(\chi_Tv_o-v_o)\|_{L^p(\R^n)}
 &\le\|\nabla v_o\|_{L^p(\R^n\setminus B_T)}\\
 &\quad+C\|v_o\|_{L^s(B_{2T}\setminus B_T)}
 \longrightarrow0
\end{aligned}
\]
as $T\to\infty$. Since $\chi_Tv_o\in C_c^\infty(\R^n)$,
we conclude that $v_o\in\D^{1,p}(\R^n)$.
\end{proof}

Apply the sharp Euclidean Sobolev inequality
to $v_o$ and use~\eqref{eq:D-stereo}.
Averaging the resulting inequality over $o$
and using Tonelli's theorem, we obtain,
for every real-valued $u\in C^\infty(\mathbb S^n)$,
\begin{equation}\label{eq:stereo-averaged}
 N(u)\le\kappa\int_{\mathbb S^n}
 \mathbb E_o
 |\nabla u(y)+\beta u(y)Z_o(y)|^p\,dV(y).
\end{equation}

\subsection{Concentrating test functions}

For $p>2$ and $a=2$, the logarithmic expansion
of the deficit for a fixed cutoff radius appears in
Frank--Jin--Wang--Xiong~\cite[Section 3, Step 5]{FJWX}.
We derive estimates that remain valid when the
cutoff radius varies with the concentration scale.

For $\theta>0$ and $u\in W^{1,p}(\mathbb S^n)$, set
\[
 \mathcal D_\theta(u)
 =\|u\|_s^\theta-K^\theta\|\nabla u\|_p^\theta.
\]
In particular, $\mathcal D_p(u)=N(u)-\kappa E(u)$.

\begin{proposition}[Estimates for cut-off Talenti functions]
\label{prop:bubbles}
Fix $n\ge2$, $1<p<n$, and $R_0\in(0,1/2)$.
Let $\chi:[0,\infty)\to[0,1]$ be a nonincreasing
Lipschitz function equal to one on $[0,1]$
and to zero on $[2,\infty)$.
For $0<\lambda<R\le R_0$, set
\[
 b_{\lambda,R}(r)=\chi(r/R)U_\lambda(r),
 \qquad \delta=\frac{\lambda}{R},
\]
where $r=d(o,\cdot)$ for a fixed pole
$o\in\mathbb S^n$.
Then $b_{\lambda,R}\in W^{1,p}(\mathbb S^n)$.

We keep $n,p,\chi,R_0$ fixed.
All asymptotic statements below are understood
in the joint limit $\lambda\downarrow0$ and
$\delta=\lambda/R\to0$, subject to any additional
conditions specified.
Unless otherwise indicated, the constants in
the estimates are independent of $R$.
\begin{align}
 a>2:\quad \mathcal D_p(b_{\lambda,R})
 &=c_{n,p}\lambda^2+O(\delta^a)+o(\lambda^2),
 \label{eq:bubble-second}\\
 a=2:\quad \mathcal D_p(b_{\lambda,R})
 &=\gamma_{n,p}\lambda^2\log(R/\lambda)
       +O(\delta^2)+O(\lambda^2),
 \label{eq:bubble-log}
\end{align}
where
\[
 c_{n,p}=\frac{(n-1)pa}{2n(a-2)}\mathsf A^{\nu-1}\mathsf A_2>0,
 \qquad \mathsf A_2=\int_{\R^n}|x|^2U^s\,dx
 \quad(a>2),
\]
and $\gamma_{n,p}=\frac{n-1}{6}\kappa\sigma2^p>0$ when $a=2$.
In particular, the $o(\lambda^2)$ term
in~\eqref{eq:bubble-second}, divided by $\lambda^2$,
tends to zero along every sequence satisfying
$\lambda\to0$ and $\delta\to0$.
For all $a>0$, we have
\[
 E(b_{\lambda,R})\to\mathsf E,\qquad
 \int_{\mathbb S^n}b_{\lambda,R}^s\,dV\to\mathsf A.
\]
The $L^p$ mass satisfies
\begin{align}
 a>p:&\quad L(b_{\lambda,R})\sim\lambda^p\int_{\R^n}U^p\,dx,
 \label{eq:bubble-mass-integrable}\\
 a=p:&\quad L(b_{\lambda,R})=
       \sigma\lambda^p\log(R/\lambda)+O(\lambda^p),
 \label{eq:bubble-mass-log}\\
 a<p,\ R\to0:&\quad L(b_{\lambda,R})\sim
       C_\chi\lambda^aR^{p-a},\qquad
 C_\chi=\sigma\int_0^2\chi(t)^pt^{p-a-1}\,dt.
 \label{eq:bubble-mass-tail}
\end{align}
If $a<p$ and $R$ is fixed, then as $\lambda\downarrow0$,
$\lambda^{-a}L(b_{\lambda,R})$ converges instead
to the finite positive limit
\[
 \sigma R^{p-a}\int_0^2\chi(t)^pt^{p-a-1}
          \left(\frac{\sin(Rt)}{Rt}\right)^{n-1}\,dt.
\]
For fixed $R$ and $a<2$, we have
$|\mathcal D_p(b_{\lambda,R})|\le C_R\lambda^a$.
More generally, for every $a>0$, fixed $R$,
and fixed $\theta>0$, we have
\begin{equation}\label{eq:bubble-coarse}
 |\mathcal D_\theta(b_{\lambda,R})|
 \le C_{R,\theta}\lambda^{\min(2,a)}(1+|\log\lambda|).
\end{equation}
\end{proposition}

\begin{proof}
Write $b=b_{\lambda,R}$ and set
$J(r)=(\sin r/r)^{n-1}$, so that the volume element
in geodesic polar coordinates is
$dV=J(r)r^{n-1}\,dr\,d\omega$.
Differentiating the radial profile of $U$, we obtain
\[
 U'(r)=-ar^{q-1}(1+r^q)^{-n/p}.
\]
Consequently,
\[
 U_\lambda(r)\le\lambda^{a/p}r^{-a},
 \qquad
 |U_\lambda'(r)|
 \le a\lambda^{a/p}r^{-a-1}.
\]
Using $n=p+(p-1)a$, we obtain
\[
 \sigma\int_R^\infty U_\lambda^s r^{n-1}\,dr
 =O(\delta^{a+q}),
 \qquad
 \sigma\int_R^\infty |U_\lambda'|^p r^{n-1}\,dr
 =O(\delta^a).
\]
By the product rule, for a.e. $r\in(R,2R)$,
\[
 |b'(r)|\le |U_\lambda'(r)|+C R^{-1}U_\lambda(r).
\]
The preceding pointwise estimates give
\[
 \int_R^{2R}|b'(r)|^p r^{n-1}\,dr=O(\delta^a).
\]
Since $b=U_\lambda$ on $[0,R]$,
$0\le b\le U_\lambda$, and $b=0$ on $[2R,\pi]$,
the preceding tail and annular estimates,
together with $0<J(r)\le1$ on $(0,2R)$, give
\[
\begin{aligned}
 \int_{\mathbb S^n}b^s\,dV
 &=\mathsf A+
 \sigma\int_0^R(J(r)-1)U_\lambda^s r^{n-1}\,dr
 +O(\delta^{a+q}),\\
 E(b)
 &=\mathsf E+
 \sigma\int_0^R(J(r)-1)|U_\lambda'|^p r^{n-1}\,dr
 +O(\delta^a).
\end{aligned}
\]
The $O$-terms in these two formulas are uniform
for $0<\lambda<R\le R_0$.

Put $c_g=(n-1)/6$. For $0<r<1$, we have
\[
 J(r)=1-c_gr^2+O_n(r^4).
\]
Suppose first that $a>2$. Then both
\[
 \mathsf A_2=\int_{\R^n}|x|^2U^s\,dx,\qquad
 \mathsf E_2=\int_{\R^n}|x|^2|\nabla U|^p\,dx
\]
are finite. After the change of variables $r=\lambda t$,
dominated convergence gives
\[
\begin{aligned}
 \int_{\mathbb S^n}b^s\,dV
 &=\mathsf A-c_g\mathsf A_2\lambda^2
   +o(\lambda^2)+O(\delta^{a+q}),\\
 E(b)
 &=\mathsf E-c_g\mathsf E_2\lambda^2
   +o(\lambda^2)+O(\delta^a).
\end{aligned}
\]
Indeed, extend $(J(\lambda t)-1)/\lambda^2$
by zero from $(0,R/\lambda)$ to $(0,\infty)$.
Since $R/\lambda\to\infty$, these extensions converge
pointwise to $-c_gt^2$ and are bounded in absolute
value by $Ct^2$.
The finiteness of $\mathsf A_2$ and $\mathsf E_2$
therefore allows us to apply dominated convergence
along any sequence in the joint limit.

Since the critical mass tends to $\mathsf A>0$,
we may apply Taylor's formula to $x^\nu$ near
$\mathsf A$. The quadratic remainder is
$O(\lambda^4+\delta^{2(a+q)})$ and hence is absorbed
in $o(\lambda^2)+O(\delta^a)$.
Using $\mathsf A^\nu=\kappa\mathsf E$, we obtain
\[
 \mathcal D_p(b)
 =c_g\bigl(\kappa\mathsf E_2
       -\nu\mathsf A^{\nu-1}\mathsf A_2\bigr)\lambda^2
       +o(\lambda^2)+O(\delta^a).
\]

For $k=0,2$, the function
$r^{n+k}(1+r^q)^{-(n-1)}$ vanishes at both endpoints,
since $a>2$. Integrating its derivative and using
$q(n-1)-n=a$, we obtain
\[
\begin{aligned}
 &\int_0^\infty r^{n+k+q-1}(1+r^q)^{-n}\,dr\\
 &\qquad=\frac{n+k}{a-k}
 \int_0^\infty r^{n+k-1}(1+r^q)^{-n}\,dr,
 \qquad k=0,2.
\end{aligned}
\]
Thus
\[
 \mathsf E=na^{p-1}\mathsf A,\qquad
 \mathsf E_2=\frac{a^p(n+2)}{a-2}\mathsf A_2.
\]
Since $\kappa\mathsf E=\mathsf A^\nu$, it follows that
\[
\begin{aligned}
 \kappa\mathsf E_2-\nu\mathsf A^{\nu-1}\mathsf A_2
 &=\mathsf A^{\nu-1}\mathsf A_2
 \left(\frac{a(n+2)}{n(a-2)}-\nu\right)\\
 &=\frac{3pa}{n(a-2)}\mathsf A^{\nu-1}\mathsf A_2.
\end{aligned}
\]
This proves~\eqref{eq:bubble-second} with the stated coefficient.

Suppose next that $a=2$.
The second moment of $U^s$ remains finite, but
that of $|\nabla U|^p$ diverges logarithmically.
More precisely,
\[
 \sigma\int_0^R r^{n+1}|U_\lambda'(r)|^p\,dr
 =\sigma2^p\lambda^2\log(R/\lambda)
 +O(\lambda^2).
\]
Indeed, after scaling, the integrand
$\sigma t^{n+1}|U'(t)|^p$ equals
$\sigma2^p/t+O(t^{-1-q})$ for $t\ge1$.
Moreover,
\[
 \int_0^R r^{n+3}|U_\lambda'(r)|^p\,dr
 \le C\lambda^4+
 C\lambda^2\int_\lambda^R r\,dr
 =O(\lambda^2).
\]
Substituting these estimates and the Taylor expansion
of $J$ into the mass and energy identities, we obtain
\[
\begin{aligned}
 \int_{\mathbb S^n}b^s\,dV
 &=\mathsf A+O(\lambda^2)+O(\delta^{2+q}),\\
 E(b)
 &=\mathsf E-c_g\sigma2^p\lambda^2\log(R/\lambda)
   +O(\lambda^2)+O(\delta^2).
\end{aligned}
\]
Applying the mean value theorem to $x^\nu$
near $\mathsf A$ and using
$\mathsf A^\nu=\kappa\mathsf E$
and $\delta^{2+q}\le\delta^2$, we obtain
\[
 \mathcal D_p(b)
 =\gamma_{n,p}\lambda^2\log(R/\lambda)
 +O(\lambda^2)+O(\delta^2),
 \qquad
 \gamma_{n,p}=c_g\kappa\sigma2^p.
\]
Finally, if $a<2$, then $|J(r)-1|\le Cr^a$
for $0<r<1$. Since
$\int_{\R^n}|x|^aU^s\,dx<\infty$, we have
\[
 \int_{\mathbb S^n}b^s\,dV
 =\mathsf A+O(\lambda^a)+O(\delta^{a+q}).
\]
For the energy, using $|J(r)-1|\le Cr^2$
and splitting the integral at $r=\lambda$, we obtain
\[
\begin{aligned}
 |E(b)-\mathsf E|
 &\le C\lambda^2
 +C\lambda^a\int_\lambda^Rr^{1-a}\,dr
 +C\delta^a\\
 &\le C(\lambda^a+\delta^a).
\end{aligned}
\]
For fixed $R$, we have $\delta^a=R^{-a}\lambda^a$.
Applying the mean value theorem to $x^\nu$ near
$\mathsf A$ and using $\mathsf A^\nu=\kappa\mathsf E$,
we obtain $|\mathcal D_p(b)|\le C_R\lambda^a$.

The preceding estimates also give
$E(b)\to\mathsf E$ and $\int b^s\,dV\to\mathsf A$
in the joint limit, since
$\lambda^2\log(R/\lambda)\to0$.

We now turn to the estimates for $L(b)$.
If $a>p$, then $U^p\in L^1(\R^n)$. After scaling,
\[
 \lambda^{-p}L(b)=\sigma\int_0^{2R/\lambda}
 \chi(\delta t)^p U(t)^pJ(\lambda t)t^{n-1}\,dt.
\]
Dominated convergence gives~\eqref{eq:bubble-mass-integrable}.

If $a=p$, the contributions to $L(b)$ from
$0<r<\lambda$ and $R<r<2R$ are both $O(\lambda^p)$,
by scaling and the pointwise bound for $U_\lambda$,
respectively.
Since $b=U_\lambda$ on $(\lambda,R)$, we have
\[
 U_\lambda(r)^p r^{n-1}
 =\lambda^p r^{-1}
   \bigl[1+O((\lambda/r)^q)\bigr].
\]
Using $J(r)=1+O(r^2)$, we obtain
\[
\begin{aligned}
 L(b)
 &=\sigma\lambda^p\log(R/\lambda)\\
 &\quad+O\left(
 \lambda^p\left[
 1+\int_\lambda^R
 \left((\lambda/r)^q+r^2\right)\frac{dr}{r}
 \right]\right)\\
 &=\sigma\lambda^p\log(R/\lambda)+O(\lambda^p).
\end{aligned}
\]

If $a<p$, the exact formula is
\[
 \frac{L(b)}{\lambda^aR^{p-a}}
 =\sigma\int_0^2\chi(t)^pt^{p-a-1}
  (1+(\delta/t)^q)^{-p\beta}
  \left(\frac{\sin(Rt)}{Rt}\right)^{n-1}\,dt.
\]
Since $a<p$, the function
$\chi(t)^pt^{p-a-1}$ is integrable on $(0,2)$.
Since $0<J(Rt)\le1$, the error from replacing
$(1+(\delta/t)^q)^{-p\beta}$ by one is bounded
in absolute value by an integral independent of $R$,
which tends to zero by dominated convergence.
Hence
\[
 \frac{L(b)}{\lambda^aR^{p-a}}
 =\sigma\int_0^2\chi(t)^pt^{p-a-1}J(Rt)\,dt
 +o(1)
\]
as $\delta\to0$, uniformly for $R\le R_0$.
The integral is bounded above and away from zero
uniformly in $R$. Letting $R\to0$ gives
\eqref{eq:bubble-mass-tail}, while keeping $R$ fixed
gives the stated limit of $\lambda^{-a}L(b_{\lambda,R})$.

For fixed $R$ and $\theta>0$, the convergences
$N(b)\to\mathsf A^\nu$ and
$\kappa E(b)\to\mathsf A^\nu$ imply, by the mean
value theorem, that
\[
 |\mathcal D_\theta(b)|
 \le C_\theta|\mathcal D_p(b)|.
\]
Combining this inequality with~\eqref{eq:bubble-second},
\eqref{eq:bubble-log}, and the fixed-radius deficit
estimate for $a<2$, we obtain~\eqref{eq:bubble-coarse}.
\end{proof}

\section{Sharp growth of the lower-order coefficient}
\label{sec:remainder}

We prove Theorem~\ref{thm:remainder} by establishing
matching lower and upper bounds for $B_\varepsilon(n,p)$
in Regions I--IV.
The lower bounds are obtained using the concentrating
test functions in Proposition~\ref{prop:bubbles},
while the upper bound in the critical case $n=p^2$
is treated separately from the other three cases.

\subsection{Lower bounds}

\begin{proposition}\label{prop:four-lower}
Let $p>2$ and $n\ge3p-2$. Then there exist constants
$c,\varepsilon_0>0$, depending only on $n,p$, such that
\[
 B_\varepsilon(n,p)\ge c f_{n,p}(\varepsilon)
 \qquad\text{for }0<\varepsilon<\varepsilon_0,
\]
where $f_{n,p}$ is defined in Theorem~\ref{thm:remainder}.
\end{proposition}

\begin{proof}
Fix $R_0=1/4$ and a cutoff function $\chi$ as in
Proposition~\ref{prop:bubbles}, chosen independently
of $n$ and $p$.
Write $b=b_{\lambda,R}$ and $\delta=\lambda/R$.
By~\eqref{eq:Beps-def},
\[
 B_\varepsilon(n,p)\ge
 \frac{\mathcal D_p(b)-\varepsilon E(b)}{L(b)}
 \qquad(0<\lambda<R\le R_0).
\]

First consider Regions I--III, where $a>2$.
In Regions I and III, take $R=R_0$, so that
$\delta^a=o(\lambda^2)$ as $\lambda\downarrow0$.
In Region II, take
\[
 R=M\lambda^{(a-2)/a},
\]
where $M>0$ is to be chosen independently of $\lambda$.
Since $\delta^a=M^{-a}\lambda^2$ and the implicit
constant in the $O(\delta^a)$ term
of~\eqref{eq:bubble-second} is independent of $R$,
we may fix $M$ sufficiently large that this term
is bounded in absolute value by $c_{n,p}\lambda^2/4$.
With these choices, $0<\lambda<R\le R_0$ for
sufficiently small $\lambda$, and $\delta\to0$
as $\lambda\downarrow0$.
Hence~\eqref{eq:bubble-second} gives
\[
 \mathcal D_p(b)\ge\frac{c_{n,p}}2\lambda^2
\]
for sufficiently small $\lambda$.
Since $E(b)\to\mathsf E$, we may now take
$\lambda=T\sqrt{\varepsilon}$ with $T$ fixed and
sufficiently large to obtain
\[
 \mathcal D_p(b)-\varepsilon E(b)
 \ge\frac{c_{n,p}}4\lambda^2
\]
for sufficiently small $\varepsilon>0$.

In Regions I, II, and III, we have $a>p$,
$2<a<p$, and $a=p$, respectively.
Moreover, $R\to0$ in Region II.
The $L^p$ mass estimates~\eqref{eq:bubble-mass-integrable},
\eqref{eq:bubble-mass-tail}, and
\eqref{eq:bubble-mass-log} thus give
\[
 L(b)\le C
 \begin{cases}
 \lambda^p,&\text{in Region I},\\
 \lambda^{2+p(a-2)/a},&\text{in Region II},\\
 \lambda^p\log(1/\lambda),&\text{in Region III}.
 \end{cases}
\]
Here we used
$\lambda^aR^{p-a}=M^{p-a}\lambda^{2+p(a-2)/a}$
in Region II.
Substituting these mass estimates and the lower bound
for $\mathcal D_p(b)-\varepsilon E(b)$ into
\eqref{eq:Beps-def}, with $\lambda=T\sqrt{\varepsilon}$,
we obtain the desired lower bounds in Regions I--III.
In particular, we use
$\log(1/\lambda)\sim\frac12\log(1/\varepsilon)$
in Region III and the identity
\[
 \frac{a-2}{a}=\frac{n-3p+2}{n-p}
\]
in Region II.

In Region IV, we have $a=2<p$. Set
\[
 \tau=\log(1/\varepsilon),\qquad
 \lambda=\sqrt{\varepsilon},\qquad
 R=M\tau^{-1/2},
\]
where $M>0$ will be chosen independently of
$\varepsilon$.
Then
\[
 \log(R/\lambda)
 =\tfrac12\tau-\tfrac12\log\tau+\log M,
 \qquad
 \delta^2=M^{-2}\varepsilon\tau.
\]
For each fixed $M>0$, as $\varepsilon\downarrow0$,
we have $R\to0$ and $\delta\to0$.
The joint asymptotic estimate~\eqref{eq:bubble-log},
together with the boundedness of $E(b)$, yields
\[
\begin{aligned}
 \mathcal D_p(b)-\varepsilon E(b)
 &=\frac{\gamma_{n,p}}2\varepsilon\tau
   +O(M^{-2}\varepsilon\tau)
   +o_M(\varepsilon\tau),
\end{aligned}
\]
where the implicit constant in the $O$-term is
independent of $M$, and the $o_M$-term is taken
as $\varepsilon\downarrow0$ with $M$ fixed.
We may therefore fix $M$ sufficiently large and
then take $\varepsilon$ sufficiently small to obtain
\[
 \mathcal D_p(b)-\varepsilon E(b)
 \ge c\varepsilon\tau.
\]
By~\eqref{eq:bubble-mass-tail},
\[
 L(b)\asymp\lambda^2R^{p-2}
 \asymp_M\varepsilon\tau^{-(p-2)/2}.
\]
Consequently,
\[
 B_\varepsilon(n,p)\ge c\tau^{p/2}
 =c[\log(1/\varepsilon)]^{p/2}.
\]
\end{proof}

\subsection{Upper bounds}

\begin{corollary}[An explicit upper bound for $p>2$]
\label{cor:drift-large-p}
Let $2<p<n$, and set
\[
 \tau_p=2^{\max\{p-3,0\}},
\]
\[
 D_2=\kappa\frac{p(p-1)\tau_p}{2}
       \beta^2\frac n{n-2},\qquad
 D_p=\kappa\tau_p\beta^pM_{n,p},
\]
where $M_{n,p}$ is defined in~\eqref{eq:stereo-moments}.
Then, for every $u\in W^{1,p}(\mathbb S^n)$,
\begin{equation}\label{eq:D-vector-global}
 N(u)\le\kappa E(u)
       +D_2E(u)^{(p-2)/p}L(u)^{2/p}
       +D_pL(u).
\end{equation}
Consequently, for every $\varepsilon>0$,
\begin{equation}\label{eq:D-Beps}
 B_\varepsilon(n,p)\le D_p+
 \frac2pD_2^{p/2}
 \left(\frac{p-2}{p}\right)^{(p-2)/2}
 \varepsilon^{-(p-2)/2}.
\end{equation}
\end{corollary}

\begin{proof}
For $X,Y\in\R^n$, the function $t\mapsto|X+tY|^p$
is of class $C^2$ since $p>2$, and its second
derivative is bounded above by
$p(p-1)|X+tY|^{p-2}|Y|^2$.
Using
\[
 |X+tY|^{p-2}
 \le\tau_p\bigl(|X|^{p-2}
       +t^{p-2}|Y|^{p-2}\bigr),
 \qquad 0\le t\le1,
\]
in Taylor's formula with integral remainder, we obtain
\[
\begin{aligned}
 |X+Y|^p
 &\le |X|^p+p|X|^{p-2}X\cdot Y\\
 &\quad+\frac{p(p-1)\tau_p}{2}
       |X|^{p-2}|Y|^2+\tau_p|Y|^p,
\end{aligned}
\]
where we used $\int_0^1(1-t)\,dt=1/2$ and
$\int_0^1(1-t)t^{p-2}\,dt=1/[p(p-1)]$.

Let $u\in C^\infty(\mathbb S^n)$ be real-valued.
Applying the preceding inequality with
$X=\nabla u(y)$ and $Y=\beta u(y)Z_o(y)$,
and using~\eqref{eq:stereo-averaged}, we obtain
\[
 N(u)\le\kappa E(u)
 +D_2\int_{\mathbb S^n}|\nabla u|^{p-2}|u|^2\,dV
 +D_pL(u).
\]
The linear term vanishes upon averaging over $o$,
since $\mathbb E_oZ_o(y)=0$.
The coefficients $D_2$ and $D_p$ follow from the
moment identities in Lemma~\ref{lem:stereographic}.
By H\"older's inequality,
\[
 \int_{\mathbb S^n}|\nabla u|^{p-2}|u|^2\,dV
 \le E(u)^{(p-2)/p}L(u)^{2/p}.
\]
This proves~\eqref{eq:D-vector-global} for smooth
functions, and the general case follows by density.

Writing $t=E(u)/L(u)$ for $u\ne0$ and
applying~\eqref{eq:D-vector-global} to the quotient
in~\eqref{eq:Beps-def}, we obtain
\[
 B_\varepsilon(n,p)
 \le D_p+\sup_{t\ge0}
 \left(D_2t^{(p-2)/p}-\varepsilon t\right).
\]
The supremum on the right is attained at
\[
 t_*=\left(\frac{(p-2)D_2}{p\varepsilon}\right)^{p/2}.
\]
Substituting this value yields~\eqref{eq:D-Beps}.
\end{proof}

\begin{proposition}[Upper bounds in Regions I, II, and IV]
\label{prop:three-upper}
The upper estimates in Theorem~\ref{thm:remainder}
hold in Regions I, II, and IV.
\end{proposition}

\begin{proof}
The upper bound in Region I follows directly
from~\eqref{eq:D-Beps}.

For Regions II and IV, Lemma~\ref{lem:weight} gives
the local inequality
\[
 N(f)\le A_RE(f)
 \qquad\text{for every }f\in W^{1,p}_0(B_R),
\]
where $A_R=\kappa(1+\omega_R)^{p(n-1)/n}$.
Since $\omega_R\to0$ as $R\downarrow0$, we have
$A_R-\kappa\le C\omega_R$ for sufficiently small $R$.
Using the bounds for $\omega_R$ in
Lemma~\ref{lem:weight}, we obtain
\[
 A_R-\kappa\le
 \begin{cases}
 CR^{2a/(a-2)},&a>2,\\
 CR^2\exp(-c_1/R^2),&a=2,
 \end{cases}
\]
for sufficiently small $R$.

Whenever $0<R\le1$ and $A_R\le\kappa+\varepsilon$,
Proposition~\ref{prop:localization} applies with
$A=\kappa+\varepsilon$ and gives
\[
 B_\varepsilon(n,p)\le
 C(\kappa+\varepsilon)R^{-p}\le CR^{-p}
 \qquad(0<\varepsilon\le1).
\]

In Region II, we have $2<a<p$. Choose
$R=c\varepsilon^{(a-2)/(2a)}$, with $c>0$
sufficiently small and fixed. Then
$A_R\le\kappa+\varepsilon$ for sufficiently small
$\varepsilon$, and hence
\[
 B_\varepsilon(n,p)\le
 C\varepsilon^{-p(a-2)/(2a)}
 =C\varepsilon^{-p(n-3p+2)/(2(n-p))}.
\]

In Region IV, we have $a=2$. Set
$\tau=\log(1/\varepsilon)$ and choose
$R^2=c_1/(2\tau)$. Then
\[
 A_R-\kappa\le C\tau^{-1}\varepsilon^2
 \le\varepsilon
\]
for sufficiently small $\varepsilon$. Consequently,
\[
 B_\varepsilon(n,p)\le CR^{-p}
 \le C[\log(1/\varepsilon)]^{p/2}.
\]
\end{proof}

\begin{proposition}[Upper bound in Region III]
\label{prop:critical-log}
Let $p>2$ and $n=p^2$. Then there exist constants
$C>0$ and $\varepsilon_0\in(0,1)$, depending only
on $n,p$, such that
\[
 B_\varepsilon(n,p)\le
 C\frac{\varepsilon^{-(p-2)/2}}{\log(1/\varepsilon)}
 \qquad(0<\varepsilon<\varepsilon_0).
\]
\end{proposition}

\begin{proof}
Suppose that the claimed upper bound fails. Then
there exists a sequence $\varepsilon_j\downarrow0$
such that, writing $B_j=B_{\varepsilon_j}(n,p)$,
\[
 B_j\varepsilon_j^{(p-2)/2}
 \log(1/\varepsilon_j)\longrightarrow+\infty.
\]
In particular, $B_j\to+\infty$.

Fix a pole $o\in\mathbb S^n$.
Applying Lemma~\ref{lem:duality} with
$A=\kappa+\varepsilon_j$, we obtain a nonnegative
extremizer $u_j\in W^{1,p}(\mathbb S^n)$ for the
inequality defining $B_j$.
Each $u_j$ is radial and nonincreasing about $o$
and satisfies $\int_{\mathbb S^n}u_j^s\,dV=1$.
Writing $E_j=E(u_j)$ and $L_j=L(u_j)$, we have
\begin{equation}\label{eq:III-extremizer}
\begin{gathered}
 (\kappa+\varepsilon_j)E_j+B_jL_j=1,\\
 -\Delta_pu_j+\frac{B_j}{\kappa+\varepsilon_j}u_j^{p-1}
 =\frac{1}{\kappa+\varepsilon_j}u_j^{s-1}
 \qquad\text{weakly on }\mathbb S^n.
\end{gathered}
\end{equation}
In particular,
\[
 E_j\le S,\qquad L_j\le B_j^{-1}\longrightarrow0.
\]

Fix $R_0=1/2$ and a radial nonincreasing Lipschitz
cutoff $0\le\eta\le1$, equal to one on $B_{R_0/2}$
and vanishing outside $B_{3R_0/4}$, with
$|\nabla\eta|\le4/R_0$. We choose $\eta$ independently
of $n,p$.

Since $B_j>0$ and $(\kappa+\varepsilon_j)^{-1}\le S$,
the normalized solutions $u_j$ satisfy the hypotheses
of the cutoff estimates in the proof of
Proposition~\ref{prop:localization}.
Hence~\eqref{eq:cutoff-two-sided} gives
\begin{equation}\label{eq:III-cutoff}
\begin{gathered}
 |E(\eta u_j)-E_j|\le CL_j,\\
 \int_{\mathbb S^n}(\eta u_j)^s\,dV=1+O(L_j).
\end{gathered}
\end{equation}

Let $v_j$ be the Euclidean transplant of $\eta u_j$
under the radial transformation in Section~\ref{sec:transplant},
with $R=R_0$. Then $v_j\in\D^{1,p}(\R^n)$ is
nonnegative, radial and nonincreasing, and vanishes
outside the fixed Euclidean ball of radius
$T_*=t(3R_0/4)$. Writing $W(x)=W_{R_0}(|x|)$ for the radial weight
on $\R^n\setminus\{0\}$, we obtain
from~\eqref{eq:transplant} and~\eqref{eq:III-cutoff}
\begin{equation}\label{eq:III-transform}
\begin{gathered}
 \int_{\R^n}|\nabla v_j|^p\,dx=E(\eta u_j),\\
 \int_{\R^n}Wv_j^s\,dx=1+O(L_j),\qquad
 \int_{\R^n}Wv_j^p\,dx\le L_j.
\end{gathered}
\end{equation}
The weight $W$ extends continuously to the origin
with $W(0)=1$ and is bounded above and away from
zero on $\{|x|\le T_*\}$.
Since $a=p>2$, the estimate
$F_{R_0}(r)=1+O(r^2)$ from the proof of
Lemma~\ref{lem:weight}, together with $t(r)\sim r$,
gives $W(x)=1+O(|x|^2)$ near the origin.

Set
\[
 A_j=\int_{\R^n}v_j^s\,dx,\qquad
 Q_j=\kappa\int_{\R^n}|\nabla v_j|^p\,dx-A_j^\nu.
\]
The quantity $Q_j$ is the Euclidean Sobolev deficit
of $v_j$ and is nonnegative by the sharp Euclidean
Sobolev inequality. Moreover,
\[
 A_j^\nu\le\kappa E(\eta u_j)
 \le\kappa(E_j+CL_j)\le1+o(1),
\]
so that $(A_j)$ is bounded.

Since $W$ is bounded away from zero on the support
of $v_j$, we have
$\int_{\R^n}v_j^p\,dx\le CL_j$.
By radial monotonicity, for each fixed $t_0>0$,
we have $v_j(t)\le C_{t_0}L_j^{1/p}$
for a.e. $t\ge t_0$. Combining this pointwise estimate with the bound
for $\int_{\R^n}v_j^p\,dx$, we obtain
\[
 \int_{|x|\ge t_0}v_j^s\,dx
 \le C_{t_0}L_j^{s/p}\longrightarrow0.
\]
Using also the boundedness of $W$ on the support
of $v_j$, we obtain
\[
 \left|A_j-\int_{\R^n}Wv_j^s\,dx\right|
 \le A_j\sup_{|x|\le t_0}|1-W(x)|
       +C_{t_0}L_j^{s/p}.
\]
Letting first $j\to\infty$ and then $t_0\downarrow0$,
we conclude from~\eqref{eq:III-transform} that
$A_j\to1$. Since
\[
 SA_j^\nu\le E(\eta u_j)
 \le E_j+CL_j\le S+o(1),
\]
we also obtain $E_j\to S$ and $Q_j\to0$.

Let $\mathcal M$ denote the family of extremizers
in~\eqref{eq:euclidean-extremals}, and put
\[
 d_j=\|\nabla v_j\|_{L^p(\R^n)},\qquad
 m_j=\|v_j\|_{L^s(\R^n)}.
\]
For all sufficiently large $j$, we have $m_j>0$.
Applying the gradient stability theorem of
Figalli--Zhang~\cite[Theorem 1.1, p.~2410]{FZ}
with stability exponent $p>2$, and then using the
sharp Euclidean Sobolev inequality for $v_j-g$,
we obtain
\[
 \inf_{g\in\mathcal M}\|v_j-g\|_{L^s(\R^n)}^p
 \le C d_j^p\left(\frac{d_j}{m_j}-K^{-1}\right).
\]
By convexity of $z\mapsto z^p$, we have
\[
 Q_j=K^pd_j^p-m_j^p
 \ge pKm_j^p\left(\frac{d_j}{m_j}-K^{-1}\right).
\]
Since $d_j$ is bounded and $m_j\to1$,
the two preceding estimates give
\[
 \inf_{g\in\mathcal M}\|v_j-g\|_{L^s(\R^n)}^p
 \le CQ_j.
\]
As $v_j$ is nonzero and compactly supported for
all sufficiently large $j$, it is not an extremizer
of the Euclidean Sobolev inequality. Thus $Q_j>0$,
and we may choose $h_j\in\mathcal M$ such that
$\|v_j-h_j\|_{L^s(\R^n)}^p\le CQ_j$.

Here $\star$ denotes symmetric decreasing
rearrangement in $\R^n$ about the origin.
Since $v_j$ is nonnegative, radial and nonincreasing,
the nonexpansivity of this rearrangement
in $L^s(\R^n)$ (see
Brock--Solynin~\cite[Theorem 3.1 and Remark 3.1]{BS2000})
gives
\[
 \|v_j-(|h_j|)^\star\|_{L^s(\R^n)}
 \le\|v_j-|h_j|\|_{L^s(\R^n)}
 \le\|v_j-h_j\|_{L^s(\R^n)}.
\]
Since $(|h_j|)^\star$ is a Talenti function centered
at the origin, there exist $c_j,\lambda_j>0$ such that
\begin{equation}\label{eq:III-FZ}
 v_j=g_j+e_j,\qquad
 g_j=(|h_j|)^\star=c_jU_{\lambda_j},\qquad
 \|e_j\|_{L^s(\R^n)}^p\le CQ_j.
\end{equation}
By~\eqref{eq:III-FZ}, we have
$\|e_j\|_{L^s(\R^n)}\to0$.
Since $\|v_j\|_{L^s(\R^n)}\to1$ and
$\|g_j\|_{L^s(\R^n)}=c_j\mathsf A^{1/s}$,
it follows that $c_j\to\mathsf A^{-1/s}$.

Since $v_j$ vanishes outside the Euclidean ball
of radius $T_*$, we have $e_j=-g_j$ almost everywhere
on $\{|x|\ge T_*\}$. By scaling,
\[
 c_j^s\int_{|x|\ge T_*/\lambda_j}U^s\,dx
 \le\|e_j\|_{L^s(\R^n)}^s\longrightarrow0.
\]
Together with $c_j\to\mathsf A^{-1/s}>0$, this
implies that $\lambda_j\to0$.

To bound the weighted critical mass from above
in terms of $A_j$, set $w=(W-1)_+$ and
$J_j=\int_{\R^n}wv_j^s\,dx$.
On the fixed ball $\{|x|\le T_*\}$ containing
the support of $v_j$, we have
\[
 0\le w(x)\le C\min\{|x|^2,1\}.
\]
Restricting the integral defining $J_j$ to the
support of $v_j$ and using the decomposition
in~\eqref{eq:III-FZ}, we obtain
\begin{equation}\label{eq:III-J}
 J_j\le C\int_{\R^n}|x|^2g_j^s\,dx
       +C\|e_j\|_{L^s(\R^n)}^s
 \le C\lambda_j^2+CQ_j^{s/p}.
\end{equation}
For the last inequality, we used the error estimate
in~\eqref{eq:III-FZ} and the scaling identity
\[
 \int_{\R^n}|x|^2g_j^s\,dx
 =c_j^s\lambda_j^2\mathsf A_2,
\]
where $\mathsf A_2=\int_{\R^n}|x|^2U^s\,dx$
is finite since $a=p>2$.
Recall also that $c_j\to\mathsf A^{-1/s}$.

By~\eqref{eq:III-transform} and the definition of $J_j$,
\[
 1-CL_j\le\int_{\R^n}Wv_j^s\,dx\le A_j+J_j.
\]
Since $A_j\to1$, the mean value theorem gives
$1-A_j^\nu\le C(L_j+J_j)$. Combining
\eqref{eq:III-extremizer}, \eqref{eq:III-cutoff},
and~\eqref{eq:III-J}, we obtain
\begin{equation}\label{eq:III-signed}
\begin{aligned}
 B_jL_j+\varepsilon_jE_j
 &=1-\kappa E_j\\
 &=1-A_j^\nu-Q_j
       +\kappa\bigl(E(\eta u_j)-E_j\bigr)\\
 &\le CL_j+C\lambda_j^2+CQ_j^{s/p}-Q_j.
\end{aligned}
\end{equation}
Since $s/p>1$ and $Q_j\to0$, we have
$CQ_j^{s/p}\le Q_j/2$ for all sufficiently large $j$.
Since $B_j\to\infty$, we also have
$CL_j\le B_jL_j/2$ for such $j$.
Absorbing these two error terms
in~\eqref{eq:III-signed} and using $E_j\to S>0$,
we obtain
\begin{equation}\label{eq:III-scale}
 B_jL_j+c\varepsilon_j+Q_j\le C\lambda_j^2.
\end{equation}
By~\eqref{eq:III-scale} and~\eqref{eq:III-FZ},
we obtain
\[
 \lambda_j\ge c\sqrt{\varepsilon_j},
 \qquad
 \|e_j\|_{L^s(\R^n)}\le C\lambda_j^{2/p}.
\]

We next compare $v_j$ with $g_j$ on a shrinking
Euclidean ball to obtain a lower bound for $L_j$.
Set $T_j=\lambda_j^{(p-2)/p}$ and
$\Omega_j=\{x\in\R^n:|x|<T_j\}$.
Since $\lambda_j\to0$, we have $T_j>\lambda_j$
for all sufficiently large $j$.
By H\"older's inequality and $1/p-1/s=1/n$,
\[
 \|e_j\|_{L^p(\Omega_j)}
 \le CT_j\|e_j\|_{L^s(\R^n)}
 \le C\lambda_j.
\]
Using the explicit formula for $U_{\lambda_j}$
with $a=p$ and $\beta=p-1$, and recalling that
$c_j\to\mathsf A^{-1/s}>0$, we obtain
$g_j(r)\ge c\lambda_jr^{-p}$ for $r\ge\lambda_j$
and all sufficiently large $j$.
Since $n=p^2$, it follows that
\[
 \|g_j\|_{L^p(\Omega_j)}^p
 \ge c\lambda_j^p
       \int_{\lambda_j}^{T_j}\frac{dr}{r}
 \ge c\lambda_j^p\log(1/\lambda_j).
\]
Here $\log(T_j/\lambda_j)=(2/p)\log(1/\lambda_j)$.

By the reverse triangle inequality and the bound
$\int_{\R^n}v_j^p\,dx\le CL_j$ established above,
\[
\begin{aligned}
 CL_j^{1/p}
 &\ge \|v_j\|_{L^p(\Omega_j)}\\
 &\ge \|g_j\|_{L^p(\Omega_j)}
       -\|e_j\|_{L^p(\Omega_j)}\\
 &\ge c\lambda_j[\log(1/\lambda_j)]^{1/p}
       -C\lambda_j.
\end{aligned}
\]
Since $\log(1/\lambda_j)\to\infty$, the last error
term is absorbed by the leading term on the right.
Consequently,
\begin{equation}\label{eq:III-tail}
 L_j\ge c\lambda_j^p\log(1/\lambda_j)
\end{equation}
for all sufficiently large $j$.

Combining~\eqref{eq:III-scale}
and~\eqref{eq:III-tail}, we obtain
\[
 B_j\le C\frac{\lambda_j^{-(p-2)}}{\log(1/\lambda_j)}.
\]
The function $x\mapsto x^{-(p-2)}/\log(1/x)$
is decreasing for sufficiently small $x>0$.
Since $\lambda_j\to0$ and
$\lambda_j\ge c\sqrt{\varepsilon_j}$,
we may use this monotonicity together with
$\log(1/(c\sqrt{\varepsilon_j}))
\sim\frac12\log(1/\varepsilon_j)$ to obtain
\[
 B_j\le
 C\frac{\varepsilon_j^{-(p-2)/2}}
        {\log(1/\varepsilon_j)}.
\]
This contradicts the choice of $(\varepsilon_j)$
and completes the proof.
\end{proof}

\begin{proof}[Proof of Theorem~\ref{thm:remainder}]
The lower bounds follow from
Proposition~\ref{prop:four-lower}.
The upper bounds are established in
Propositions~\ref{prop:three-upper}
and~\ref{prop:critical-log}.
\end{proof}

\section{Optimal powers}
\label{sec:C-outer}

We now prove Theorem~\ref{thm:powers}.
We first establish the required inequalities for
$1<p<2$ and $n>p^2$ by combining local inequalities
with height truncation.
We then derive necessary bounds on the admissible
powers using concentrating test functions.
For $p>2$ and $n\ge3p-2$, we combine the known
failure of the sharp $p$-power inequality with
Theorem~\ref{thm:remainder} and Lemma~\ref{lem:power}
to complete the classification.

\subsection{Local inequalities and height truncation}

\begin{lemma}[Local sharp Sobolev inequalities]
\label{lem:C-local}
For each fixed pair $(n,p)$, there exists $R\in(0,1)$
with the following property: for every $\theta$
in the appropriate range below,
\[
 \begin{cases}
 \theta>0,&a<2,\\
 0<\theta\le\min\{2,s\},&a>2,\\
 0<\theta<\min\{2,s\},&a=2.
 \end{cases}
\]
there is a finite constant $C_\theta\ge0$ such that
\[
 \|f\|_s^\theta
 \le K^\theta\|\nabla f\|_p^\theta
      +C_\theta\|f\|_p^\theta
 \qquad\text{for all }f\in W^{1,p}_0(B_R).
\]
If $a<2$, one may take $C_\theta=0$.
\end{lemma}

\begin{proof}
If $a<2$, the sharp Dirichlet Sobolev inequality
on a sufficiently small ball follows from
Bandle--Peletier--Stingelin~\cite[Theorem 2.1
and pp.~1389--1390]{BPS}.
Raising it to any positive power gives the result.

Suppose $a\ge2$ and fix $R=1/2$.
We may assume that $f\ne0$.
Let $w=f^\star$ and let $v$ be its Euclidean
transplant as in Section~\ref{sec:transplant}.
The proof of Lemma~\ref{lem:weight} gives
$(F_R(r)-1)_+\le Cr^2$ if $a>2$, and
$(F_R(r)-1)_+\le Cr^2\log(R/r)$ if $a=2$.
Since $F_R$ is bounded and
$W_R(t(r))=F_R(r)^{(n+a)/a}$, it follows that
$(W_R(t(r))-1)_+\le C_\theta r^\theta$
for every $\theta$ in the stated range.
Using $(W_R-1)_+/W_R\le(W_R-1)_+$,
the mass and energy identities
in~\eqref{eq:transplant}, and the sharp Euclidean
Sobolev inequality, we obtain
\begin{equation}\label{eq:C-weighted}
\begin{aligned}
 \|w\|_s^s
 &\le \|v\|_{L^s(\R^n)}^s+
 \int_{B_R}
 \frac{(W_R(t(r))-1)_+}{W_R(t(r))}w^s\,dV\\
 &\le K^s\|\nabla w\|_p^s
       +C_\theta\int_{B_R}r^\theta w^s\,dV.
\end{aligned}
\end{equation}

By radial monotonicity and $|B_r|\ge cr^n$,
we have $r^nw(r)^p\le C\|w\|_p^p$ for a.e.
$r\in(0,R)$. Since $n(s-p)/p=s$, it follows that
$r^sw(r)^{s-p}\le C\|w\|_p^{s-p}$.
Multiplying by $w^p$ and integrating gives
$\int_{B_R}r^sw^s\,dV\le C\|w\|_p^s$.
Applying H\"older's inequality to
$r^\theta w^s=(r^sw^s)^{\theta/s}(w^s)^{1-\theta/s}$
for $0<\theta<s$, and including the preceding
estimate at $\theta=s$, we obtain
\begin{equation}\label{eq:C-moment}
 \int_{B_R}r^\theta w^s\,dV
 \le C_\theta\|w\|_p^\theta\|w\|_s^{s-\theta},
 \qquad 0<\theta\le s.
\end{equation}

If $K\|\nabla w\|_p\ge\|w\|_s$, the desired
inequality is immediate.
Otherwise, using $x^s\le x^\theta$ for
$0\le x\le1$ and $0<\theta\le s$, together
with~\eqref{eq:C-weighted} and~\eqref{eq:C-moment},
we obtain
\[
\begin{aligned}
 \|w\|_s^\theta-K^\theta\|\nabla w\|_p^\theta
 &\le
 \frac{\|w\|_s^s-K^s\|\nabla w\|_p^s}
      {\|w\|_s^{s-\theta}}\\
 &\le C_\theta\|w\|_p^\theta.
\end{aligned}
\]
The result for $f$ follows by rearrangement.
\end{proof}

\begin{lemma}[Localization by height truncation]
\label{lem:C-height}
Let $R\in(0,\pi)$ and $0<\theta\le\min\{s,n/p\}$.
Suppose that there exists $C_{\rm loc}\ge0$ such that
\[
 \|f\|_s^\theta
 \le K^\theta\|\nabla f\|_p^\theta
      +C_{\rm loc}\|f\|_p^\theta
 \qquad\text{for all }f\in W^{1,p}_0(B_R).
\]
Then there exists a finite constant $C>0$ such that
\[
 \|u\|_s^\theta
 \le K^\theta\|\nabla u\|_p^\theta
      +C\|u\|_p^\theta
 \qquad\text{for all }u\in W^{1,p}(\mathbb S^n).
\]
No restriction $\theta\le p$ is imposed.
\end{lemma}

\begin{proof}
By homogeneity and spherical decreasing rearrangement,
it suffices to consider a nonnegative function
$w\in W^{1,p}(\mathbb S^n)$, radial and nonincreasing
about the center of $B_R$, with $\|w\|_s=1$.
Set $z=\|w\|_p$.

Using the locally absolutely continuous representative
of the radial profile on $(0,\pi)$, set
$h=w(R)$ and $v=(w-h)_+$.
By the chain rule for Sobolev functions,
$v\in W^{1,p}(\mathbb S^n)$.
Radial monotonicity implies that $v=0$ almost everywhere
outside $B_R$, so its restriction belongs to
$W^{1,p}_0(B_R)$.
Moreover,
\[
 h\le |B_R|^{-1/p}z,\qquad
 \|v\|_p\le z,\qquad
 \|\nabla v\|_p\le\|\nabla w\|_p.
\]

Since $0\le w^s-v^s\le shw^{s-1}$, we have
$1-\|v\|_s^s\le Cz\int_{\mathbb S^n}w^{s-1}\,dV$.
By H\"older's inequality, we obtain
\[
 \int_{\mathbb S^n}w^{s-1}\,dV
 \le
 \begin{cases}
 Cz^{s-1},&s-1\le p,\\
 z^{p/(s-p)},&s-1>p,
 \end{cases}
\]
where the first estimate uses the finite volume of
$\mathbb S^n$, and the second follows by interpolation
between $L^p$ and $L^s$, using $\|w\|_s=1$.
Since $1+p/(s-p)=n/p$ and
$s-1\le p$ if and only if $s\le n/p$, we obtain
\begin{equation}\label{eq:C-height-loss}
 0\le1-\|v\|_s^s\le Cz^k,
 \qquad k=\min\{s,n/p\}.
\end{equation}

As $\theta\le k\le s$, $\|v\|_s\le1$,
and $z\le V_n^{1/n}$, it follows that
$1-\|v\|_s^\theta\le1-\|v\|_s^s
\le C_\theta z^\theta$.
Applying the local inequality to $v$, we obtain
\[
 1\le\|v\|_s^\theta+C_\theta z^\theta
 \le K^\theta\|\nabla w\|_p^\theta
      +(C_{\rm loc}+C_\theta)z^\theta.
\]
The result follows by homogeneity and rearrangement.
\end{proof}

\begin{proposition}[Admissible powers for $1<p<2$]
\label{prop:C-low-p}
Let $1<p<2$. If $n\le p^2$, then the sharp
$p$-power inequality holds. If $n>p^2$, then the
sharp inequality holds at the power
\[
 \theta_*=\min\{2,s,n/p\}>p.
\]
In both cases, all smaller positive powers
are also admissible.
\end{proposition}

\begin{proof}
For $n\le p^2$, the sharp $p$-power inequality
follows from Aubin~\cite[Theorem 8, pp.~588--591]{Aubin}.

Suppose $n>p^2$. Each of $2,s,n/p$ is strictly
larger than $p$, so $\theta_*>p$.
If $a=2$, then $n=3p-2<2p$, whence
$n/p<2$ and $s>n\ge2$.
Thus $\theta_*<\min\{2,s\}$ in this case.
Thus, in each of the cases $a<2$, $a=2$, and $a>2$,
Lemma~\ref{lem:C-local} yields the local inequality
at $\theta_*$.
Since $\theta_*\le\min\{s,n/p\}$, applying
Lemma~\ref{lem:C-height} to this local inequality
gives the global inequality at $\theta_*$.

In both cases $n\le p^2$ and $n>p^2$,
the admissibility of all smaller positive powers
follows by subadditivity.
\end{proof}

\subsection{Necessary bounds from concentrating tests}

Fix a cutoff function $\chi$ and a radius
$R\in(0,1/2)$ as in Proposition~\ref{prop:bubbles},
and write $b_\lambda=b_{\lambda,R}$.
Unless otherwise specified, we keep $n,p,\chi,R$
fixed as $\lambda\downarrow0$.
Since $b_\lambda(r)=U_\lambda(r)$ for $0\le r\le\lambda<R$
and $\sin r\asymp r$ in this region, scaling and
the identity $\beta s=n$ give
\begin{equation}\label{eq:C-coupling}
 \int_{\mathbb S^n}b_\lambda^{s-1}\,dV
 \ge c\lambda^\beta.
\end{equation}

For any family $(u_\lambda)$ in $W^{1,p}(\mathbb S^n)$
such that both $N(u_\lambda)$ and $\kappa E(u_\lambda)$
converge to $\mathsf A^\nu$ as $\lambda\downarrow0$,
the mean value theorem applied to $x\mapsto x^{\theta/p}$
gives, for each fixed $\theta>0$,
\begin{equation}\label{eq:C-convert}
 \mathcal D_\theta(u_\lambda)
 =\left[
 \frac{\theta}{p}(\mathsf A^\nu)^{\theta/p-1}+o(1)
 \right]\mathcal D_p(u_\lambda).
\end{equation}

\begin{proposition}[Tests with a constant background]
\label{prop:C-background}
Every $\theta\in\mathcal T(n,p)$ satisfies
the following bounds:
\begin{enumerate}
 \item $\theta\le s$;
 \item $\theta\le n/p$ if $n>p^2$;
 \item $\theta\le p$ if $a<2$ and $a\le p$.
\end{enumerate}
\end{proposition}

\begin{proof}
Fix $\theta\in\mathcal T(n,p)$.
By the fixed-radius estimates in
Proposition~\ref{prop:bubbles}, we have
$\|b_\lambda\|_p\to0$ and
$\mathcal D_p(b_\lambda)\ge-C\lambda^a$.
The latter bound follows from the positivity of
$\mathcal D_p(b_\lambda)$ for sufficiently small
$\lambda$ if $a\ge2$, and from the fixed-radius
deficit estimate if $a<2$.
Applying~\eqref{eq:C-convert} with the outer power $s$,
we obtain $\mathcal D_s(b_\lambda)\ge-C\lambda^a$.

Adding a constant does not change the gradient.
For all sufficiently small $\lambda$, both
$\|b_\lambda+h\|_s$ and $K\|\nabla b_\lambda\|_p$
lie in a fixed compact subinterval of $(0,\infty)$,
uniformly for $h\in[0,1]$.
By the admissibility of $\theta$ and the mean value
theorem applied to $x\mapsto x^{s/\theta}$, we obtain
\[
 \mathcal D_s(b_\lambda+h)
 \le C_\theta(\|b_\lambda\|_p+h)^\theta.
\]
On the other hand, the gradient terms cancel in
$\mathcal D_s(b_\lambda+h)-\mathcal D_s(b_\lambda)$.
Using the pointwise inequalities
$(b+h)^s-b^s\ge h^s$ and
$(b+h)^s-b^s\ge shb^{s-1}$ for $b,h\ge0$,
together with~\eqref{eq:C-coupling}, we obtain
\[
 \begin{aligned}
 c(h^s+h\lambda^\beta)
 &\le \mathcal D_s(b_\lambda+h)
       -\mathcal D_s(b_\lambda)\\
 &\le C\lambda^a+
       C_\theta(\|b_\lambda\|_p+h)^\theta,
 \end{aligned}
\]
where the constants are independent of $h$ and
$\lambda$ for $h\in[0,1]$ and sufficiently small
$\lambda>0$.

Fixing $h\in(0,1]$ and letting $\lambda\downarrow0$
in the preceding estimate, we obtain
$ch^s\le C_\theta h^\theta$.
Letting $h\downarrow0$ gives $\theta\le s$.

If $n>p^2$, then $a>p$ and
$\|b_\lambda\|_p\le C\lambda$
by~\eqref{eq:bubble-mass-integrable}.
Taking $h=\lambda$ and using $1+\beta=n/p$, we obtain
$\lambda^{n/p}\le C_\theta(\lambda^a+\lambda^\theta)$.
Since $a-n/p=(n-p^2)/(p(p-1))>0$,
this implies $\theta\le n/p$.

Finally, suppose that $a<2$ and $a\le p$.
Take $h=T\lambda^{a/p}$, where $T>0$ is to be fixed.
Since $a/p+\beta=a$ and the constants in the preceding
estimate are independent of $h$, we may fix $T$
sufficiently large to absorb the term $C\lambda^a$.
For sufficiently small $\lambda$, we then have
$h\le1$ and obtain
\[
 \lambda^a\le
 C_\theta(\|b_\lambda\|_p+T\lambda^{a/p})^\theta.
\]
The fixed-radius mass estimates give
$\|b_\lambda\|_p\le C\lambda^{a/p}$ if $a<p$,
and
$\|b_\lambda\|_p\le
C\lambda[\log(1/\lambda)]^{1/p}$ if $a=p$.
In either case, for sufficiently small $\lambda$,
we have
$\lambda^a\le C_{\theta,T}
\lambda^{a\theta/p}[\log(1/\lambda)]^{\theta/p}$.
If $\theta>p$, the right-hand side is $o(\lambda^a)$,
a contradiction.
Thus $\theta\le p$.
\end{proof}

\begin{proposition}[Curvature tests]
\label{prop:C-curvature}
If $a>2$ and $n>p^2$, then $\theta\le2$
for every $\theta\in\mathcal T(n,p)$.
If $a=2$ and $p<2$, then on every fixed small
geodesic ball $B_\rho$ there is no finite constant $C$
such that
\[
 \|f\|_s^2\le K^2\|\nabla f\|_p^2+C\|f\|_p^2
 \qquad\text{for all }f\in W^{1,p}_0(B_\rho).
\]
\end{proposition}

\begin{proof}
For the second assertion, take the pole to be
the center of $B_\rho$ and choose the fixed
cutoff radius $R$ so that
$0<R<\min\{1/2,\rho/2\}$.
Then $b_\lambda$ is compactly supported in $B_\rho$.
In both cases $a>p$, so
$\|b_\lambda\|_p\asymp\lambda$
by~\eqref{eq:bubble-mass-integrable}.
Combining~\eqref{eq:bubble-second},
\eqref{eq:bubble-log} and~\eqref{eq:C-convert}
gives
\[
 \frac{\mathcal D_\theta(b_\lambda)}
      {\|b_\lambda\|_p^\theta}
 \ge c_\theta
 \begin{cases}
  \lambda^{2-\theta},&a>2,\\
  \lambda^{2-\theta}\log(1/\lambda),&a=2,
 \end{cases}
\]
for every fixed $\theta>0$ and sufficiently small
$\lambda$. If $a>2$, the quotient diverges for
$\theta>2$, proving the first assertion.
For the second assertion, taking $\theta=2$
gives a logarithmically divergent quotient
as $\lambda\downarrow0$.
Since $b_\lambda\in W^{1,p}_0(B_\rho)$,
this excludes any finite constant $C$
in the local inequality.
\end{proof}

\subsection{Proof of Theorem~\ref{thm:powers}}

\begin{proposition}\label{prop:C-p-two}
For every integer $n\ge3$, one has
$\mathcal T(n,2)=(0,2]$ and $B_0(n,2)=V_n^{-2/n}$.
\end{proposition}

\begin{proof}
The sharp inequality~\eqref{eq:quadratic-intro}
and testing with constant functions give
$B_0(n,2)=V_n^{-2/n}$.
By subadditivity, $(0,2]\subset\mathcal T(n,2)$.

To prove the reverse inclusion, fix $o\in\mathbb S^n$
and consider the conformal extremizers
\[
 u_\lambda(y)=\mathsf A^{-1/s}
 \left(
 \frac{\lambda}{
 1-o\cdot y+\lambda^2(1+o\cdot y)}
 \right)^\beta,
 \qquad 0<\lambda<1.
\]
The stereographic transformation in
Lemma~\ref{lem:stereographic} maps $u_\lambda$
to the Euclidean extremizer $\mathsf A^{-1/s}U_\lambda$.

By the conformal equivalence of the Euclidean and
spherical Sobolev inequalities in the quadratic case
(see Hebey~\cite[Theorem 5.1, pp.~121--124]{HebeyBook}
and Beckner~\cite[(31), (33), pp.~231--232]{Beckner1993}),
$u_\lambda$ attains equality in~\eqref{eq:quadratic-intro}.
Moreover, \eqref{eq:D-stereo} gives $\|u_\lambda\|_s=1$.
Since $s>2$, this normalization implies that
the family $(u_\lambda^2)$ is uniformly integrable.
As $u_\lambda\to0$ pointwise away from $o$
when $\lambda\downarrow0$, we conclude that
$z_\lambda:=\|u_\lambda\|_2\to0$.

Since $u_\lambda$ attains equality
in~\eqref{eq:quadratic-intro} and $\|u_\lambda\|_s=1$,
we have $K\|\nabla u_\lambda\|_2\le1$.
Thus, for every $\theta>2$,
$\mathcal D_\theta(u_\lambda)\ge
\mathcal D_2(u_\lambda)$, and hence
\[
 \frac{\mathcal D_\theta(u_\lambda)}{z_\lambda^\theta}
 \ge
 \frac{\mathcal D_2(u_\lambda)}{z_\lambda^\theta}
 =V_n^{-2/n}z_\lambda^{2-\theta}
 \longrightarrow+\infty.
\]
This rules out every $\theta>2$ and completes
the proof.
\end{proof}

\begin{proof}[Proof of Theorem~\ref{thm:powers}]

Suppose first that $1<p<2$.
If $n\le p^2$, Proposition~\ref{prop:C-low-p}
gives admissibility at $p$.
Since $a\le p<2$ in this case,
Proposition~\ref{prop:C-background} excludes
every $\theta>p$.
Hence $\mathcal T(n,p)=(0,p]$.

If $n>p^2$, Proposition~\ref{prop:C-low-p}
gives admissibility at
$\theta_*=\min\{2,s,n/p\}$.
By Proposition~\ref{prop:C-background},
every $\theta\in\mathcal T(n,p)$ satisfies
$\theta\le\min\{s,n/p\}$.
If $\theta_*=2$, then $n\ge2p>3p-2$, so $a>2$.
Proposition~\ref{prop:C-curvature} therefore gives
the additional bound $\theta\le2$.
It follows that no power greater than $\theta_*$
is admissible, and hence
$\mathcal T(n,p)=(0,\theta_*]$.

The case $p=2$ follows from
Proposition~\ref{prop:C-p-two}.

Let $p>2$. If $n<3p-2$, then
$p\in\mathcal T(n,p)$ by~\eqref{eq:known-finiteness}.
Since $a<2<p$, Proposition~\ref{prop:C-background}
excludes every $\theta>p$.
Thus $\mathcal T(n,p)=(0,p]$.

Suppose now that $n\ge3p-2$.
By~\eqref{eq:known-finiteness},
$B_0(n,p)=+\infty$.
Since admissibility at any power $\theta\ge p$
would imply admissibility at $p$, we have
$\mathcal T(n,p)\subset(0,p)$.
Combining Theorem~\ref{thm:remainder}
with Lemma~\ref{lem:power}, we obtain,
for every $0<\theta<p$,
\[
 \theta\in\mathcal T(n,p)
 \quad\Longleftrightarrow\quad
 f_{n,p}(\varepsilon)
 =O\left(\varepsilon^{-(p/\theta-1)}\right)
 \quad\text{as }\varepsilon\downarrow0.
\]

If $n\ge p^2$, this condition is equivalent to
$p/\theta-1\ge(p-2)/2$.
Indeed, at $n=p^2$, the factor $1/\log(1/\varepsilon)$
does not compensate for any strictly negative
power of $\varepsilon$.
Thus $\mathcal T(n,p)=(0,2]$.

If $3p-2<n<p^2$, the condition becomes
\[
 \theta\le
 \frac{p}{1+\frac{p(n-3p+2)}{2(n-p)}}
 =\frac{2p(n-p)}{-3p^2+np+2n}.
\]
This gives the claimed closed interval.

Finally, if $n=3p-2$, the logarithmic growth
of $B_\varepsilon(n,p)$ satisfies the preceding
growth condition for every $0<\theta<p$.
Hence $\mathcal T(n,p)=(0,p)$.

\end{proof}

\paragraph{Some endpoint cases.}
The condition $a=2$ with $1<p<2$ occurs only
for $(n,p)=(2,4/3)$ and $(3,5/3)$.
The largest admissible global powers are $3/2$
and $9/5$, respectively, even though the local
sharp inequality with outer power $2$ fails by
Proposition~\ref{prop:C-curvature}.
The quadratic endpoint must therefore be excluded
at these two parameter pairs in
Aubin--Li~\cite[Proposition 5.1, p.~370]{AL}.
Replacing the unit round metric $g$ by $k^{-1}g$
multiplies the quadratic deficit quotient by $k$,
so the same obstruction persists for arbitrarily
small positive sectional curvature $k$.
This is consistent with the strict restriction
$p<(n+2)/3$ in Remark 5.1 and Subcase 2.1
of their proof.

For $(n,p)=(3,6/5)$, we have $s=2$.
Admissibility at $\theta_*=s=2$ follows from
the weighted moment estimate~\eqref{eq:C-moment}
and the height truncation in Lemma~\ref{lem:C-height}.

\section{Estimates for the second Sobolev constant}\label{sec:D-second}

We first establish uniform estimates on spherical
caps and derive a lower bound for $B_0(n,p)$
using transplanted Talenti functions.
Together with quantitative localization, these
estimates yield parts (b) and (c) of
Theorem~\ref{thm:second}.
Part (a) follows from the quadratic case and the
explicit upper bound in Appendix~\ref{app:bounds}.

\subsection{Uniform estimates on caps}

Suppose that $0<a<2$ and $d:=p-a>0$.
For $0<R<\pi$, let $G_R$ be as
in~\eqref{eq:transplant-def} and define
\begin{equation}\label{eq:cap-data}
\begin{aligned}
 h_a(R)&=(\sin R)^{-a}
       -a\int_0^R(1-\cos\tau)(\sin\tau)^{-a-1}\,d\tau,\\
 N_R&=\sigma\int_0^RG_R(r)^p\sin^{n-1}r\,dr.
\end{aligned}
\end{equation}

The integral representation of $H_R$ in the proof of
Lemma~\ref{lem:weight} shows that
$H_R(r)\to h_a(R)$ as $r\downarrow0$.
Since $H_R$ is nondecreasing, $h_a(R)\ge0$
implies $H_R\ge0$ and hence $W_R\le1$.
Combining the transplantation identities
in~\eqref{eq:transplant} with the sharp Euclidean
Sobolev inequality and spherical rearrangement,
we obtain
\begin{equation}\label{eq:cap-sharp}
 h_a(R)\ge0
 \quad\Longrightarrow\quad
 N(f)\le\kappa E(f)
 \quad\text{for all }f\in W^{1,p}_0(B_R).
\end{equation}
Since $G_R(r)\sim r^{-a}$ as $r\downarrow0$,
we have $G_R(r)^p\sin^{n-1}r\sim r^{d-1}$,
where $d=n-ap>0$.
Moreover, $G_R(r)\to0$ as $r\uparrow R$.
Thus $0<N_R<\infty$.

\begin{lemma}[Uniform small-radius asymptotics]
\label{lem:D-beta}
Fix $n\ge4$ and $I\Subset(1,n)$.
For $p\in I$ satisfying
$3/2\le a=(n-p)/(p-1)<2$ and $d=p-a>0$,
set $\eta=2-a$.
There exists $R_0=R_0(n,I)>0$ such that,
for every $0<R<R_0$,
\begin{align}
 h_a(R)&=R^{-a}\left[1-\frac{a(a+1)}{6\eta}R^2
                         +O_{n,I}(R^4)\right],\label{eq:D-h}\\
 N_R&=\frac{\sigma}{a}R^d
       \mathrm B(d/a,p+1)\left[1+O_{n,I}(R^2)\right].\label{eq:D-NR}
\end{align}
The error bounds are uniform in $p$ under the
stated assumptions. In particular, the relative
error in~\eqref{eq:D-NR} remains uniform as $d\downarrow0$.
\end{lemma}

\begin{proof}
Taylor expansion at zero gives, uniformly for
$3/2\le a<2$,
\[
 (\sin R)^{-a}
 =R^{-a}\left(1+\frac a6R^2+O(R^4)\right),
 \qquad
 a(1-\cos r)(\sin r)^{-a-1}
 =\frac a2r^{1-a}+O(r^{3-a}).
\]
Substituting these expansions into~\eqref{eq:cap-data}
gives~\eqref{eq:D-h}, since
$a/6-a/(2\eta)=-a(a+1)/(6\eta)$ and
$\int_0^Rr^{3-a}\,dr=O(R^{4-a})$
uniformly for $3/2\le a<2$.

For $0<\tau<R$, we have
$(\sin\tau/\tau)^{-a-1}=1+O(R^2)$ uniformly
in $\tau$ and $a$.
By the positivity of the kernel in the definition
of $G_R$, it follows that
\[
 G_R(r)=(r^{-a}-R^{-a})(1+O(R^2)),
\]
uniformly in $r$ and $a$.
Since $p\in I$ and
$(\sin r/r)^{n-1}=1+O(R^2)$, we obtain
\[
 \begin{aligned}
 N_R
 &=\sigma(1+O_{n,I}(R^2))
   \int_0^R(r^{-a}-R^{-a})^p r^{n-1}\,dr\\
 &=\frac{\sigma}{a}R^d\mathrm B(d/a,p+1)
   (1+O_{n,I}(R^2)).
 \end{aligned}
\]
The beta integral evaluation in the last line
follows from the substitution $z=(r/R)^a$
and requires only $a>0$ and $d=n-ap>0$.
The relative error in the integrand is bounded
uniformly by $C(n,I)R^2$, so the same relative
error bound holds after integration, even as
$d\downarrow0$.
\end{proof}

\begin{lemma}[A lower bound from transplanted Talenti functions]
\label{lem:D-transplant}
Let $2<p<n$ and $a<2$. For every $0<R<\pi$
such that $h_a(R)<0$, one has
\begin{equation}\label{eq:D-C}
 B_0(n,p)\ge\mathcal C_{n,p}(R):=
 \nu\mathsf A^{\nu-1}\frac{\sigma}{a}
 \frac{-h_a(R)}{N_R},
\end{equation}
where $\mathsf A$ is given by~\eqref{eq:talenti-data}.
\end{lemma}

\begin{proof}
Fix $p$ and $R$. By~\eqref{eq:transplant-def},
the weight in~\eqref{eq:transplant} satisfies
\[
 W_R(t)=
 \bigl(1+t^aH_R(r(t))\bigr)^{-(n+a)/a}.
\]
Since $H_R(r)\to h_a(R)$ as $r\downarrow0$, we have
\[
 W_R(t)=1-\frac{n+a}{a}h_a(R)t^a+o(t^a)
 \qquad\text{as }t\downarrow0.
\]
By the definition of $W_R$, we have
$W_R(t)\le t^{-(n+a)}$ for all $t>0$.
Together with the preceding expansion, this gives
the global bound
$|W_R(t)-1|\le C_{p,R}t^a$ for all $t>0$.

Define $u_\lambda(r)=U_\lambda(t(r))$ for $0<r<R$
and extend $u_\lambda$ by zero to $\mathbb S^n$.
By the radial transplantation established in
Section~\ref{sec:transplant}, the restriction of
$u_\lambda$ to $B_R$ belongs to $W^{1,p}_0(B_R)$,
and $E(u_\lambda)=\mathsf E$.
Since $n+a=q(n-1)$, we have
\[
 \int_{\R^n}|x|^aU^s\,dx
 =\frac{\sigma}{q}\mathrm B(n-1,1)
 =\frac{\sigma}{n+a}.
\]
By~\eqref{eq:transplant} and scaling, we have
\[
 \frac{\int_{\mathbb S^n}u_\lambda^s\,dV-\mathsf A}
      {\lambda^a}
 =\int_{\R^n}U^s
   \frac{W_R(\lambda|x|)-1}{\lambda^a}\,dx.
\]
The integrand is bounded in absolute value by
the integrable function $C_{p,R}|x|^aU^s$.
Hence, by dominated convergence,
\[
\begin{aligned}
 \lim_{\lambda\downarrow0}
 \frac{\int_{\mathbb S^n}u_\lambda^s\,dV
       -\mathsf A}{\lambda^a}
 &= -\frac{n+a}{a}h_a(R)
      \int_{\R^n}|x|^aU^s\,dx\\
 &= -\frac{\sigma}{a}h_a(R).
\end{aligned}
\]

For $0<r<R$, we also have
$\lambda^{-a/p}u_\lambda(r)
=(\lambda^q+t(r)^q)^{-\beta}\le G_R(r)$.
As $\lambda\downarrow0$,
$\lambda^{-a/p}u_\lambda(r)\to G_R(r)$.
Since $G_R^p$ is integrable on $B_R$,
dominated convergence gives
$\lambda^{-a}L(u_\lambda)\to N_R$.
Since $E(u_\lambda)=\mathsf E$ and
$\kappa\mathsf E=\mathsf A^\nu$,
the two mass limits and the mean value theorem
for $x\mapsto x^\nu$ give,
by~\eqref{eq:Beps-def} with $\varepsilon=0$,
\[
 B_0(n,p)\ge
 \lim_{\lambda\downarrow0}
 \frac{N(u_\lambda)-\kappa E(u_\lambda)}
      {L(u_\lambda)}
 =\nu\mathsf A^{\nu-1}\frac{\sigma}{a}
       \frac{-h_a(R)}{N_R}.
\]
All limits are taken with $p$ and $R$ fixed;
the resulting bound may subsequently be applied
with $R=R(p)$.
\end{proof}

\subsection{Parameter limits}

\begin{proposition}[Estimates near the critical parameters]
\label{prop:D-critical}
The following assertions hold.
\begin{enumerate}[label=(\roman*)]
\item For each fixed $n\ge5$, put $P=(n+2)/3$. Then
\begin{equation}\label{eq:D-match}
 B_0(n,p)\asymp_n(p-P)^{-P/2}
 \qquad\text{as }p\downarrow P.
\end{equation}
Moreover, define
\[
 q_c=\frac P{P-1},\qquad
 \nu_c=\frac{n-P}{n},\qquad
 \mathsf A_c=\frac{\sigma}{q_c}\mathrm B(n/q_c,n/P)
\]
and
\begin{equation}\label{eq:D-Ln}
 \mathcal L_n=
 \frac{\nu_c\mathsf A_c^{\nu_c-1}}
      {\mathrm B((P-2)/2,P+1)}
 \left(\frac{P-1}{3}\right)^{P/2}
 \frac2{P-2}\left(\frac{P-2}{P}\right)^{P/2}.
\end{equation}
Then
\begin{equation}\label{eq:D-liminf}
 \liminf_{p\downarrow P}(p-P)^{P/2}B_0(n,p)
 \ge\mathcal L_n>0.
\end{equation}

\item In dimension four,
\begin{equation}\label{eq:D-four}
 \liminf_{p\downarrow2}B_0(4,p)
 \ge\frac43V_4^{-1/2}>B_0(4,2).
\end{equation}
\end{enumerate}
\end{proposition}

\begin{proof}
Fix $n\ge4$ and put $P=(n+2)/3$.
For $p>P$, set $\delta=p-P$ and $\eta=2-a$.
Then
\[
 \eta=\frac{3\delta}{p-1},
 \qquad a\to2,\qquad d\to P-2
 \quad\text{as }p\downarrow P.
\]
Choose a compact neighborhood $I$ of $P$
in $(1,n)$.
For $p>P$ sufficiently close to $P$,
we have $p\in I$, $3/2\le a<2$, and $d>0$.

Combining~\eqref{eq:D-h}, \eqref{eq:D-NR}
and~\eqref{eq:D-C}, and using $a+d=p$, we obtain
\[
 \mathcal C_{n,p}(R)=
 \frac{\nu\mathsf A^{\nu-1}}{\mathrm B(d/a,p+1)}
 R^{-p}
 \frac{\dfrac{a(a+1)}{6\eta}R^2-1+O_{n,I}(R^4)}
      {1+O_{n,I}(R^2)}
\]
for $0<R<R_0(n,I)$ with $h_a(R)<0$,
where $R_0(n,I)$ may be taken sufficiently small.
The error bounds are uniform in $p$,
including as $d\downarrow0$.

Suppose first that $n\ge5$.
Taking $R_+=\frac12\sqrt\eta$ in~\eqref{eq:D-h}
gives $R_+^ah_a(R_+)\to3/4$ as $p\downarrow P$.
Thus $h_a(R_+)>0$ for $p$ sufficiently close to $P$,
and~\eqref{eq:cap-sharp} gives the sharp local
Sobolev inequality on $B_{R_+}$.
By the uniform localization estimate in
Proposition~\ref{prop:localization}, we obtain
\[
 B_0(n,p)\le C_nR_+^{-p}
 \asymp_n\delta^{-P/2},
\]
where we used $\delta^{-(p-P)/2}\to1$.

For the lower bound, fix $c>1$ and take
$R=c\sqrt\eta$.
By~\eqref{eq:D-h}, we have
$R^ah_a(R)\to1-c^2<0$, so
$h_a(R)<0$ for $p$ sufficiently close to $P$.
Moreover,
\[
 \delta^{P/2}R^{-p}
 \longrightarrow
 c^{-P}\left(\frac{P-1}{3}\right)^{P/2}.
\]
As $p\downarrow P$, the quantities $q,\nu,\mathsf A$
converge to $q_c,\nu_c,\mathsf A_c$, respectively,
and $d/a\to(P-2)/2>0$.
Using the preceding asymptotic formula for
$\mathcal C_{n,p}(R)$ and the lower bound
in Lemma~\ref{lem:D-transplant}, we obtain
\[
 \liminf_{p\downarrow P}\delta^{P/2}B_0(n,p)
 \ge
 \frac{\nu_c\mathsf A_c^{\nu_c-1}}
      {\mathrm B((P-2)/2,P+1)}
 \left(\frac{P-1}{3}\right)^{P/2}
 (c^2-1)c^{-P}.
\]
Choosing $c^2=P/(P-2)$, which maximizes
$(c^2-1)c^{-P}$ for $c>1$, proves~\eqref{eq:D-liminf}.
Combining this lower bound with the upper bound
above yields~\eqref{eq:D-match}.

Suppose now that $n=4$, so that $P=2$ and
$d/\eta=(p+2)/3\to4/3$.
Take $R=\eta^{1/4}$.
Since $R^2/\eta\to\infty$, it follows
from~\eqref{eq:D-h} that $h_a(R)<0$
for $p$ sufficiently close to $2$.
As $d/a\to0$, the beta--gamma identity and the
recurrence formula for the gamma function give
\[
 \frac da\mathrm B(d/a,p+1)
 =\frac{\Gamma(1+d/a)\Gamma(p+1)}
        {\Gamma(p+1+d/a)}
 \longrightarrow1.
\]
Using this limit together with $a(a+1)/6\to1$
and $R^2/\eta\to\infty$ in the preceding formula
for $\mathcal C_{4,p}(R)$, we obtain
\[
 \mathcal C_{4,p}(R)
 \sim\frac{\nu}{a}\mathsf A^{\nu-1}
       \frac d\eta R^{2-p}.
\]
By~\eqref{eq:talenti-data},
$\mathsf A\to(\sigma/2)\mathrm B(2,2)=V_4/16$.
Since $\nu\to1/2$ and $a\to2$, we have
$(\nu/a)\mathsf A^{\nu-1}\to V_4^{-1/2}$.
Moreover, since $\eta\sim3(p-2)$, we obtain
$R^{2-p}=\exp(-(p-2)\log\eta/4)\to1$.
It follows that
$\mathcal C_{4,p}(R)\to\frac43V_4^{-1/2}$.
Lemma~\ref{lem:D-transplant} therefore gives
the asserted lower bound for the liminf.
Since $B_0(4,2)=V_4^{-1/2}$ by
Proposition~\ref{prop:C-p-two},
this proves~\eqref{eq:D-four}.
\end{proof}

\begin{proof}[Proof of Theorem~\ref{thm:second}]
For part (a), Proposition~\ref{prop:C-p-two} gives
$B_0(n,2)=V_n^{-2/n}$ for $n\ge3$, while
Proposition~\ref{prop:D-poles} gives the
left-hand asymptotic estimate and the formula
for the left derivative.
Parts (b) and (c) follow from
Proposition~\ref{prop:D-critical}.
In particular, the strict inequality
in~\eqref{eq:D-four} implies discontinuity
at $p=2$ in dimension four.
\end{proof}

\appendix
\section{Effective bounds}\label{app:bounds}

We collect explicit bounds for the second Sobolev constant.
For $p\le2$, the upper bound yields the left-hand
asymptotic estimates in Theorem~\ref{thm:second}.
For $p>2$, we derive explicit lower bounds by
perturbing constant functions and, when $n<3p-2$,
obtain upper and lower bounds using geodesic balls.

\subsection{Bounds for \texorpdfstring{$p\le2$}{p <= 2}}

\begin{proposition}[An explicit upper bound and left-hand asymptotics]
\label{prop:D-poles}
Let $n\geq2$, $1<p<n$, and $p\leq2$.  Then
\begin{equation}\label{eq:D-polebound}
 V_n^{-p/n}\leq B_0(n,p)\leq\mathcal U_{n,p}:=
 \kappa\left(\frac{n-p}{p}\right)^p
 \frac{\Gamma((n+p)/2)\Gamma((n-p)/2)}{\Gamma(n/2)^2}.
\end{equation}
Consequently, for each fixed $n\geq3$, one has
\begin{equation}\label{eq:D-left}
 B_0(n,p)=V_n^{-p/n}+O_n((2-p)^2)\quad(p\uparrow2),\qquad
 B'_{0,-}(n,2)=-\frac{\log V_n}{n}V_n^{-2/n}.
\end{equation}
\end{proposition}

\begin{proof}

For $X,Y\in\R^n$ and $1<p\le2$, the concavity and
subadditivity of $t\mapsto t^{p/2}$ give
\[
 \frac{|X+Y|^p+|X-Y|^p}{2}
 \le (|X|^2+|Y|^2)^{p/2}
 \le |X|^p+|Y|^p.
\]
For each fixed $y\in\mathbb S^n$, the distribution
of $Z_o(y)$ is symmetric about the origin, as shown
in the proof of Lemma~\ref{lem:stereographic}.
For a smooth real-valued function $u$, apply the
preceding inequality with $X=\nabla u(y)$ and
$Y=\beta u(y)Z_o(y)$ and average over $o$.
Using this symmetry and the moment identity
in~\eqref{eq:stereo-moments}, we deduce
from~\eqref{eq:stereo-averaged} that
\[
 N(u)\le\kappa E(u)+\kappa\beta^pM_{n,p}L(u).
\]
Since $\kappa\beta^pM_{n,p}=\mathcal U_{n,p}$,
density gives the upper bound in~\eqref{eq:D-polebound}.
The lower bound follows by testing with constant functions.

To study the behavior of $B_0(n,p)$ as $p\uparrow2$,
fix $n\ge3$.
The explicit expression for $\mathcal U_{n,p}$
in~\eqref{eq:D-polebound} defines a positive analytic
function of $p$ in a neighborhood of $2$.
Set
\begin{equation}\label{eq:pole-log-comparison}
 F(p)=\log(\mathcal U_{n,p}V_n^{p/n}),\qquad
 D(p)=\log\frac{\Gamma(n+1)}
 {\sigma\Gamma(n/p)\Gamma(n+1-n/p)}.
\end{equation}
Writing $m(p)=\log M_{n,p}$, we obtain
from~\eqref{eq:euclidean-constant}
\[
 \begin{aligned}
 F(p)={}&\log\frac{n-p}{n}
 +(p-1)\log(p-1)-p\log p\\
 &+\frac pn[D(p)+\log V_n]+m(p).
 \end{aligned}
\]
Comparing~\eqref{eq:euclidean-constant} at $p=2$
with~\eqref{eq:quadratic-intro}, we obtain
$D(2)+\log V_n=n\log2$.
Moreover, since $M_{n,2}=n/(n-2)$ by
Lemma~\ref{lem:stereographic}, we have
$\mathcal U_{n,2}=V_n^{-2/n}$ and hence $F(2)=0$.

Let $\psi$ and $\psi_1$ denote the digamma and
trigamma functions, respectively.
Differentiating $D$ and $m$ and using the recurrences
$\psi(x+1)-\psi(x)=1/x$ and
$\psi_1(x+1)-\psi_1(x)=-x^{-2}$, we obtain
\[
 \begin{aligned}
 D'(2)&=-\frac12,&
 D''(2)&=\frac34-\frac{n^2}{8}\psi_1(n/2),\\
 m'(2)&=\frac1n+\frac1{n-2},&
 m''(2)&=\frac12\psi_1(n/2)-\frac1{n^2}
              +\frac1{(n-2)^2}.
 \end{aligned}
\]
Substituting these derivatives into the expression
for $F$, we find that $F'(2)=0$ and
\[
 \chi_n:=F''(2)
 =\frac{(n+2)(n-1)}{2n^2}
 -\frac{n-2}{4}\psi_1(n/2).
\]
Since $F(2)=F'(2)=0$, Taylor expansion of the
analytic function $e^{F(p)}$ at $p=2$ gives
\begin{equation}\label{eq:pole-quadratic-expansion}
 V_n^{p/n}\mathcal U_{n,p}
 =1+\frac{\chi_n}{2}(p-2)^2+O_n(|p-2|^3).
\end{equation}
Combining~\eqref{eq:pole-quadratic-expansion}
with the two-sided bounds in~\eqref{eq:D-polebound}
gives the first estimate in~\eqref{eq:D-left}.
Since $B_0(n,2)=V_n^{-2/n}$, the first estimate
in~\eqref{eq:D-left} also yields the stated
left derivative.
\end{proof}

\begin{corollary}[An upper bound for the normalized quadratic remainder]
\label{cor:left-quadratic}
Let $n\ge3$, and let $\chi_n$ be as in the proof of
Proposition~\ref{prop:D-poles}. Then
\[
 0\le\liminf_{p\uparrow2}
 \frac{V_n^{p/n}B_0(n,p)-1}{(2-p)^2}
 \le\limsup_{p\uparrow2}
 \frac{V_n^{p/n}B_0(n,p)-1}{(2-p)^2}
 \le\frac{\chi_n}{2}.
\]
\end{corollary}

\begin{proof}
This follows directly from~\eqref{eq:D-polebound}
and~\eqref{eq:pole-quadratic-expansion}.
\end{proof}

\subsection{Perturbations of constants}

For $2<p<n$, the strict inequality
$B_0(n,p)>V_n^{-p/n}$ follows from
Hebey~\cite[Proposition 4.1, pp.~83--84]{HebeyBook};
see also Druet--Hebey~\cite[Theorem 8.11]{DH}.
We now obtain explicit lower bounds by testing
with perturbations of constant functions.
For a fixed $o\in\mathbb S^n$, define
\begin{equation}\label{eq:constant-moments}
 P_z(t)=\frac1{V_n}\int_{\mathbb S^n}(1+t\,o\cdot y)^z\,dV(y),
 \qquad z>0,\quad 0\le t\le\frac12.
\end{equation}

\begin{proposition}[An explicit improvement of the volume lower bound]
\label{prop:D-gap}
Let $n\geq3$, $2<p<n$, and set
\[
 g=\frac{p(s-p)}{n+1},\qquad
 j=\frac\sigma{V_n}\mathrm B\left(\frac{n+p}{2},\frac12\right),
 \qquad\gamma=\frac{g}{9\,3^s},
\]
\[
 t_0=\min\left\{\frac12,
       \left(\frac{\gamma}{\kappa V_n^{p/n}j}\right)^{1/(p-2)}
       \right\}.
\]
Then
\begin{equation}\label{eq:D-gap}
 B_0(n,p)\geq V_n^{-p/n}(1+\gamma t_0^2)>V_n^{-p/n}.
\end{equation}
In particular this is an effective estimate throughout the finite
region $n<3p-2$.
\end{proposition}

\begin{proof}
Let $\mu$ be the normalized volume measure
$d\mu=dV/V_n$ on $\mathbb S^n$, and write
$\mathbb E$ for expectation with respect to $\mu$.
For the function $X(y)=o\cdot y$, we have
$\mathbb EX=0$, $\Var_\mu X=1/(n+1)$ and
$\mathbb E|\nabla X|^p=j$,
where the last identity follows by integration
in polar coordinates.

For $0<t\le1/2$, set $Y=\log(1+tX)$ and
$f(z)=\log\mathbb E e^{zY}$.
Define the exponentially tilted probability measure
$d\mu_z=e^{zY-f(z)}d\mu$.
Since $1-t\le e^Y\le1+t$, we have, for $0\le z\le s$,
\[
 \frac{d\mu_z}{d\mu}
 \ge\left(\frac{1-t}{1+t}\right)^s.
\]
The derivative of $x\mapsto\log(1+tx)$ on $[-1,1]$
is bounded below by $t/(1+t)$.
Hence the mean value theorem and the independent-copy
formula for variance give
\[
 \Var_\mu Y\ge(t/(1+t))^2\Var_\mu X.
\]
Since $Y$ is bounded, differentiation under the
integral is justified.
Using the variational characterization
$\Var_{\mu_z}Y=\inf_{c\in\R}\int(Y-c)^2\,d\mu_z$
and the preceding two bounds, we obtain, for $0\le z\le s$,
\[
 f''(z)=\Var_{\mu_z}Y
 \ge\frac{t^2}{n+1}
       \frac{(1-t)^s}{(1+t)^{s+2}}
 =:m_t.
\]
The function $z\mapsto f(z)-m_tz^2/2$
is therefore convex on $[0,s]$ and vanishes at zero.
Using $p=\nu s$, we obtain
\[
 \nu f(s)-f(p)\ge\frac{m_tp(s-p)}2.
\]
Since $P_z(t)=e^{f(z)}$ for $z>0$,
exponentiating and using $e^x\ge1+x$ give
\[
 \frac{P_s(t)^\nu}{P_p(t)}
 \ge1+\frac g2t^2
       \frac{(1-t)^s}{(1+t)^{s+2}}.
\]
By~\eqref{eq:Beps-def} with $\varepsilon=0$,
$B_0(n,p)$ is bounded below by the quotient
$\mathcal T_{n,p}(t)$ obtained from $u=1+tX$.
Since $P_p(t)\ge1$ by Jensen's inequality, we have
\begin{equation}\label{eq:D-T}
\begin{aligned}
 \mathcal T_{n,p}(t)
 &:=
 \frac{V_n^{-p/n}P_s(t)^\nu-\kappa jt^p}{P_p(t)}\\
 &\ge V_n^{-p/n}
 \left(1+\frac g2t^2
       \frac{(1-t)^s}{(1+t)^{s+2}}\right)
 -\kappa jt^p.
\end{aligned}
\end{equation}
For $0<t\le1/2$, we have
$(1-t)^s/(1+t)^{s+2}\ge4/(9\,3^s)$.
By the definition of $t_0$,
$\kappa jt_0^{p-2}\le V_n^{-p/n}\gamma$.
Taking $t=t_0$ in~\eqref{eq:D-T}
therefore gives~\eqref{eq:D-gap}.
\end{proof}

\begin{proposition}[Another explicit lower bound]
\label{prop:D-oldgap}
Let $2<p<n$, and let $g,j$ be as in
Proposition~\ref{prop:D-gap}. Set
\[
 t_*=\min\left\{
 \frac1{4(s+1)},
 \left(\frac g{8\kappa V_n^{p/n}j}\right)^{1/(p-2)}
 \right\}.
\]
Then
\[
 B_0(n,p)\ge V_n^{-p/n}
 \left(1+\frac{gt_*^2}{8}\right).
\]
\end{proposition}

\begin{proof}
For $0\le t\le1/[4(s+1)]$, both $1-st$ and
$1-(s+2)t$ are positive.
Bernoulli's inequality therefore gives
\[
 \frac{(1-t)^s}{(1+t)^{s+2}}
 \ge(1-st)(1-(s+2)t)
 \ge1-2(s+1)t
 \ge\frac12.
\]
By definition, $t_*$ lies in this range and satisfies
$\kappa jt_*^{p-2}\le V_n^{-p/n}g/8$.
Substituting $t=t_*$ into~\eqref{eq:D-T}
gives the asserted lower bound.
\end{proof}

\subsection{Bounds from caps}

\begin{proposition}[Explicit cap estimates in the finite region]
\label{prop:D-effective-cap}
Let $2<p<n$ and $a<2$, so that $d=p-a>0$.
Here $m_R$ and $C_p$ are as in
Proposition~\ref{prop:localization}, while
$h_a$ and $\mathcal C_{n,p}$ are given by
\eqref{eq:cap-data} and~\eqref{eq:D-C}, respectively.
Set $c_0=6/5$ and
\[
 R_+=\min\left\{1,\sqrt{\frac{2-a}{a c_0^{a+1}}}\right\}.
\]
Then the following upper bound holds:
\begin{equation}\label{eq:D-upper}
 B_0(n,p)\leq\overline B_{n,p}:=
       (2^p-1)m_{R_+}+\kappa C_p4^pR_+^{-p}.
\end{equation}
The radius $R_-=\pi-3^{-1/a}$ satisfies $h_a(R_-)<-1$,
and hence
$B_0(n,p)\ge\mathcal C_{n,p}(R_-)$.
Set
\[
 R_b=\sqrt{6c_0^a(2-a)/a}.
\]
If $R_b\le1$, then we also have the explicit lower bound
\begin{equation}\label{eq:D-Lb}
 B_0(n,p)\geq\mathcal L_b:=
 \frac{\nu\mathsf A^{\nu-1}c_0^{a-p(a+1)}}
       {\mathrm B((p-a)/a,p+1)}R_b^{-p}.
\end{equation}
\end{proposition}

\begin{proof}
For $0<r\le1$, we have $r/c_0\le\sin r\le r$
and $r^2/3\le1-\cos r\le r^2/2$.
Using these inequalities in~\eqref{eq:cap-data},
we obtain, for $0<R\le1$,
\[
 \begin{aligned}
 R^{-a}\left[1-\frac{ac_0^{a+1}}{2(2-a)}R^2\right]
 &\le h_a(R)\\
 &\le R^{-a}\left[c_0^a-\frac{aR^2}{3(2-a)}\right].
 \end{aligned}
\]
By the choice of $R_+$, we have
$h_a(R_+)\ge R_+^{-a}/2>0$.
Thus~\eqref{eq:cap-sharp} gives the sharp local
Sobolev inequality on $B_{R_+}$.
Applying the explicit localization
estimate~\eqref{eq:explicit-upper} with $R=R_+$
gives~\eqref{eq:D-upper}.

If $R_b\le1$, the upper estimate for $h_a$ gives
$-h_a(R_b)\ge c_0^aR_b^{-a}$.
Using $\sin\tau\ge\tau/c_0$ in the definition
of $G_R$ in~\eqref{eq:transplant-def}, we obtain,
for $0<r<R\le1$,
\[
 G_R(r)\le c_0^{a+1}(r^{-a}-R^{-a}).
\]
Since $\sin r\le r$ and $p-a=n-ap>0$,
the beta integral evaluation in the proof of
Lemma~\ref{lem:D-beta} yields
\[
 N_R\le\frac{\sigma}{a}c_0^{p(a+1)}R^{p-a}
       \mathrm B((p-a)/a,p+1).
\]
Since $h_a(R_b)<0$, substituting the bounds for
$-h_a(R_b)$ and $N_{R_b}$ into~\eqref{eq:D-C}
gives~\eqref{eq:D-Lb}.

Finally, the definition~\eqref{eq:cap-data} gives
$h'_a(R)=-a(\sin R)^{-a-1}$ and
$h_a(\pi/2)\le1$.
Integrating the derivative and using
$\sin r\le\pi-r$, we obtain, for $0<t<\pi/2$,
\[
 h_a(\pi-t)\le1-t^{-a}+(\pi/2)^{-a}.
\]
Taking $t=3^{-1/a}\in(0,\pi/2)$, we obtain
$h_a(R_-)\le-2+(\pi/2)^{-a}<-1$.
Lemma~\ref{lem:D-transplant} therefore yields
$B_0(n,p)\ge\mathcal C_{n,p}(R_-)$.
\end{proof}

\begin{corollary}[Explicit two-sided bounds]
\label{cor:bracket}
For $2<p<n$ and $a<2$, use the notation of
Propositions~\ref{prop:D-gap}, \ref{prop:D-oldgap}
and~\ref{prop:D-effective-cap}, and set
\[
\begin{aligned}
 \underline B_{n,p}=\max\biggl\{&
 V_n^{-p/n}(1+\gamma t_0^2),\quad
 V_n^{-p/n}\left(1+\frac{gt_*^2}{8}\right),\\
 &\mathcal T_{n,p}(1/2),\quad
 \mathcal C_{n,p}(R_-),\quad
 \mathcal L_b\biggr\},
\end{aligned}
\]
where $\mathcal L_b$ is omitted if $R_b>1$.
Then
\[
 \underline B_{n,p}\le B_0(n,p)\le\overline B_{n,p}.
\]
Moreover, for each fixed $n\ge5$ and $P=(n+2)/3$,
\[
 \limsup_{p\downarrow P}
 \frac{\overline B_{n,p}}{\underline B_{n,p}}
 <\infty.
\]
\end{corollary}

\begin{proof}
The upper bound is~\eqref{eq:D-upper}.
Each term in the definition of $\underline B_{n,p}$
is a lower bound for $B_0(n,p)$ by
Propositions~\ref{prop:D-gap}, \ref{prop:D-oldgap},
\ref{prop:D-effective-cap} or by the test
quotient~\eqref{eq:D-T} at $t=1/2$.

Fix $n\ge5$. As $p\downarrow P$, we have
$a\to2$, $p-a\to P-2>0$ and $R_b,R_+\to0$.
In particular, $R_b\le1$ for $p$ sufficiently
close to $P$, so $\mathcal L_b$ is included
in $\underline B_{n,p}$.
Since $m_R=|B_{R/4}|^{-p/n}$ and
$|B_{R/4}|\ge c_nR^n$ for $0<R\le1$,
we have $m_R\le C_nR^{-p}$ uniformly for $p$ near $P$.
The coefficients in~\eqref{eq:D-upper}
and~\eqref{eq:D-Lb} have finite positive limits
as $p\downarrow P$. Hence
\[
 \overline B_{n,p}\le C_nR_+^{-p},
 \qquad
 \underline B_{n,p}\ge\mathcal L_b\ge c_nR_b^{-p}.
\]
For $p$ sufficiently close to $P$, the minimum
in the definition of $R_+$ is attained by the
square-root term. Hence
$R_b/R_+=\sqrt{6c_0^{2a+1}}$,
and the ratio $\overline B_{n,p}/\underline B_{n,p}$
remains bounded as $p\downarrow P$.
\end{proof}

\section*{Declarations}
\addcontentsline{toc}{section}{Declarations}

\noindent\textbf{Conflict of interest.}
The authors declare that they have no conflicts of interest.

\medskip

\noindent\textbf{Data availability.}
Data sharing is not applicable to this article, as no datasets were generated or analyzed during the current study.

\medskip

\noindent\textbf{Generative AI disclosure.}
The authors initiated this research and led it throughout its development. ChatGPT was used to assist in discussing certain mathematical details, checking mathematical arguments and references, and preparing the manuscript. The authors have carefully reviewed the AI-assisted material, further refined the exposition, and take full responsibility for the mathematical results and the final manuscript.

\phantomsection

\end{document}